\documentstyle[11pt]{article}
\title{The Deligne-Simpson problem for zero index of rigidity
\footnote{Research partially supported by INTAS grant 97-1644}}
\author{Vladimir Petrov Kostov\\ \\ 
\hspace{7cm}{\sl To the memory of my mother}} 
\date{}
\newtheorem{tm}{Theorem}
\newtheorem{lm}[tm]{Lemma}
\newtheorem{cor}[tm]{Corollary}
\newtheorem{prop}[tm]{Proposition}
\newtheorem{defi}[tm]{Definition}
\newtheorem{nota}[tm]{Notation}
\newtheorem{rem}[tm]{Remark}
\newtheorem{ex}[tm]{Example}
\begin{document}
\maketitle 

\begin{abstract}
We consider the {\em Deligne-Simpson problem}: 
{\em Give necessary and sufficient conditions for the choice of the conjugacy 
classes $c_j\subset gl(n,{\bf C})$ or $C_j\subset GL(n,{\bf C})$, 
$j=1,\ldots ,p+1$, so that there 
exist irreducible $(p+1)$-tuples of matrices $A_j\in c_j$ whose sum is 0 
or of matrices $M_j\in C_j$ whose product is $I$.}
The matrices $A_j$ (resp. $M_j$) are interepreted as matrices-residua of 
Fuchsian linear systems (resp. as monodromy operators of regular systems) on 
Riemann's sphere. 

We consider the case when the sum of the dimensions of the conjugacy classes 
$c_j$ or $C_j$ is $2n^2$ and we prove a theorem of non-existence of 
such irreducible $(p+1)$-tuples. 
\end{abstract}
\tableofcontents 

\section{Introduction}

In the present paper we sonsider a particular case of the {\em Deligne-Simpson 
problem (DSP)}:

{\em Give necessary and sufficient conditions for the choice of the conjugacy 
classes $c_j\subset gl(n,{\bf C})$ or $C_j\subset GL(n,{\bf C})$, 
$j=1,\ldots ,p+1$, so that there 
exist irreducible $(p+1)$-tuples of matrices $A_j\in c_j$ or $M_j\in C_j$ 
satisfying respectively the equality} 

\begin{equation}\label{A_j}
A_1+\ldots +A_{p+1}=0
\end{equation}
or 

\begin{equation}\label{M_j}
M_1\ldots M_{p+1}=I~.
\end{equation}
``Irreducible'' means ``not having a common proper invariant subspace'', i.e. 
impossible to conjugate simultaneously the $(p+1)$ matrices to a 
block upper-triangular form. 
The problem is connected with the theory of linear regular systems of 
differential equations on Riemann's sphere: 

\begin{equation}\label{regular} 
\dot{X}=A(t)X
\end{equation}
Here the $n\times n$-matrix $A(t)$ is meromorphic on ${\bf C}P^1$, with 
poles at the points $a_1$, $\ldots$, $a_{p+1}$; the unknown variables $X$ 
form also a matrix $n\times n$. Such a system is called {\em regular} at the 
pole $a_j$ if one has $||X(t-a_j)||=O(|t-a_j|^{N_j})$ for some 
$N_j\in {\bf R}$ when the solution is restricted to a sector of sufficiently 
small radius and centered at $a_j$.  

A particular case of a regular system is a {\em Fuchsian} one, i.e. with 
logarithmic poles:

\begin{equation}\label{Fuchs}
{\rm d}X/{\rm d}t=(\sum _{j=1}^{p+1}A_j/(t-a_j))X
\end{equation}
where $A_j\in gl(n,{\bf C})$ are its {\em matrices-residua}; in the absence 
of a pole at $\infty$ one has (\ref{A_j}). 

As a result of a linear change of variables 

\begin{equation}\label{W}
X\mapsto W(t)X
\end{equation}
the matrix $A(t)$ of a regular system (\ref{regular}) undergoes the 
{\em gauge transformation} 

\begin{equation}\label{gauge}
A(t)\mapsto -W^{-1}\dot{W}+W^{-1}A(t)W
\end{equation}
The $n\times n$-matrix $W$ is meromorphic on ${\bf C}P^1$, its poles if any 
are usually among the points $a_j$, and outside them 
det$W\not\equiv 0$. The only invariant of a regular system under the 
linear changes (\ref{W}) is its {\em monodromy group}. This is the group 
generated by the {\em monodromy operators}. 

A {\em monodromy operator} is a linear operator mapping 
the solution space of a regular system onto itself. It is defined as follows: 
one fixes a base 
point $a\neq a_j$ for $j=1,\ldots ,p+1$, the value at $a$ of the solution 
$X$, i.e. a matrix $B\in GL(n,{\bf C})$ and a closed contour $\Gamma$ 
passing through 
$a$. The monodromy operator $M$ defined by the homotopy equivalence 
class of the contour $\Gamma$ 
maps the solution $X$ 
with $X|_{t=a}=B$ onto 
the value at $a$ of its analytic continuation along the contour 
(notation: $X\stackrel{\Gamma}{\mapsto}XM$). 

Fix $(p+1)$ contours whose homotopy equivalence classes generate 
$\pi _1({\bf C}P^1\backslash \{ a_1,\ldots ,a_{p+1}\})$. One usually 
chooses the contours such that $\Gamma _j$ consists of a segment $[a,x_j]$ 
($x_j$ is close to $a_j$), 
of a small circumference (centered at $a_j$, passing through $x_j$,  
circumventing $a_j$ counterclockwise and not containing inside any other pole 
$a_i$) and 
of the segment $[x_j,a]$. We assume that for $i\neq j$ one has 
$\Gamma _i\cap \Gamma _j=\{ a\}$ and that the index of the contour increases 
when one turns around $a$ clockwise. For such a 
choice of the contours the monodromy operators 
$M_j$ satisfy the condition (\ref{M_j}). This means that one can choose as 
generators of the monodromy group any $p$ out of the $p+1$ operators $M_j$. 

The monodromy group is an 
antirepresentation of $\pi _1({\bf C}P^1\backslash \{ a_1,\ldots ,a_{p+1}\})$ 
into $GL(n,{\bf C})$ because one has $X\stackrel{\Gamma _i\Gamma _j}
{\mapsto}XM_jM_i$ (although we often write ``representation'' instead). The 
change of $a$ and $B$ changes the monodromy group 
to a conjugate one. 

If the contours defining the operators $M_j$ are 
chosen like above, then $M_j$ is conjugate to the corresponding 
{\em operator of local 
monodromy} defined by a small lace circumventing the pole $a_j$ 
counterclockwise. 
Therefore in the case of matrices $M_j$ the DSP admits the interpretation: 

{\em For which $(p+1)$-tuples of local monodromies do there exist 
irreducible monodromy groups with such local monodromies ?}

\begin{rem}\label{exp}
{\em The eigenvalues $\lambda _{k,j}$ of the matrix-residuum $A_j$ of a 
Fuchsian system are 
connected with $\sigma _{k,j}$, the ones of the monodromy 
operator $M_j$ by $\exp (2\pi i\lambda _{k,j})=\sigma _{k,j}$.}
\end{rem}

\section{Definitions and known facts}

\subsection{The quantities $d_j$, $r_j$ and $\kappa$; the construction $\Psi$; 
(poly)multiplicity vectors\protect\label{drkappa}}

\begin{defi}
{\em A Jordan normal form (JNF) of size $n$ 
is a collection of positive integers indexed by two indices -- 
$J^n=\{ b_{i,k}\}$ -- where $k$ is the index of an eigenvalue, $i$ is the 
index of the Jordan block of size $b_{i,k}$ 
with this eigenvalue; $k=1,\ldots ,\rho$, 
$i=1,\ldots ,s_k$. We assume that all $\rho$ eigenvalues are distinct and that 
for each $k$ one has $b_{1,k}\geq \ldots \geq b_{s_k,k}$.}
\end{defi}

{\bf Convention.} All Jordan matrices and Jordan blocks are presumed to be 
upper-triangular.

\begin{defi}
{\em Denote by $J(X)$ the JNF of the matrix $X$. 
We say that the DSP is {\em solvable} (resp. {\em weakly solvable}) 
for a given 
$\{ J_j^n\}$ and given eigenvalues if there exists an 
irreducible $(p+1)$-tuple (resp. a $(p+1)$-tuple with a trivial 
centralizer) of 
matrices $M_j$ satisfying (\ref{M_j}) or of matrices $A_j$ satisfying 
(\ref{A_j}), with $J(M_j)=J_j^n$ or $J(A_j)=J_j^n$ and with the given 
eigenvalues. By definition, the DSP is solvable for 
$n=1$.}
\end{defi}

For a given conjugacy class $C$ (in $gl(n,{\bf C})$ or $GL(n,{\bf C})$) 
we denote by $d(C)$ its dimension (which is 
always even) and by $r(C)$ the quantity 
min$_{\lambda \in {\bf C}}$~rk$(X-\lambda I)$ for $X\in C$. The quantity 
$n-r(C)$ is the greatest number of Jordan blocks with one and the same 
eigenvalue. We set $d_j=d(c_j)$ (resp. $d_j=d(C_j)$) and $r_j=r(c_j)$ 
(resp. $r_j=r(C_j)$). The quantities $r(C)$ and $d(C)$ depend not on the 
conjugacy class $C$ but only on the JNF defined by it.

The following two conditions are necessary for the existence of irreducible 
$(p+1)$-tuples of matrices $M_j$ satisfying (\ref{M_j}) or of matrices $A_j$ 
satisfying (\ref{A_j}), see \cite{Si} and \cite{Ko3}, \cite{Ko4}:

\[ \begin{array}{lll}d_1+\ldots +d_{p+1}\geq 2n^2-2&\hspace{30mm}&
(\alpha _n)\\
{\rm for~all~}j~~r_1+\ldots +\hat{r}_j+\ldots +r_{p+1}\geq n&
\hspace{30mm}&(\beta _n)\end{array}\]

\begin{defi}\label{rigindex}
{\em The quantity $\kappa =2n^2-d_1-\ldots -d_{p+1}$ is called the 
{\em index of rigidity}. If condition $(\alpha _n)$ holds, then it 
takes the values 2, 0,$-2$, $\ldots$. 
Call {\em rigid} the case $\kappa =2$ 
(i.e. for which condition $(\alpha _n)$ is an equality).}
\end{defi} 

The rigid case has been studied in \cite{Ka}. In the present paper we study 
the case $\kappa =0$. These two cases are of particular interest because 
they seem to contain all non-trivial examples when the DSP is not 
weakly solvable. (An example is called non-trivial if the JNFs $J_j^n$ 
satisfy the conditions of Theorem~\ref{generic} below.)

\begin{defi}
{\em Denote by $J_j^n$ the JNF 
of size $n$ defined by the class $c_j$ or $C_j$ and by $\{ J_j^n\}$
the $(p+1)$-tuple of these JNFs. For $n>1$ define the map  
$\Psi : \{ J_j^n\} \mapsto \{ J_j^{n_1}\}$ 
if the condition $(\beta _n)$ holds and the condition} 

\[ r_1+\ldots +r_{p+1}\geq 2n\hspace{4cm}(\omega _n)\]
{\em does not hold. Namely, set $n_1=(\sum _{j=1}^{p+1}r_j)-n$; 
hence, $n_1<n$. 
For each $j$ the new JNF $J_j^{n_1}$ is defined after $J_j^n$ by choosing an 
eigenvalue with the maximal possible number $n-r_j$ of Jordan blocks, by 
decreasing by 1 the sizes of the smallest $n-n_1$ of them and by deleting the 
Jordan blocks of size 0. One has $n-n_1\leq n-r_j$ because $(\beta _n)$ holds. 
If there are several eigenvalues with maximal number of Jordan blocks, then 
we choose any of them.}
\end{defi}

\begin{defi}
{\em A {\em multiplicity vector (MV)} is a vector whose components are 
non-negative 
integers whose sum is $n$. Notation: 
$\Lambda _j^n=(m_{1,j},\ldots ,m_{i_j,j})$, 
$m_{1,j}\geq \ldots \geq m_{i_j,j}$, $m_{1,j}+\ldots +m_{i_j,j}=n$. 
The components have the meaning of the 
multiplicities of the eigenvalues of a matrix $A_j$ or $M_j$ (for the sake of 
convenience we admit components equal to 0). A 
{\em polymultiplicity vector (PMV)} is the $(p+1)$-tuple of MVs defined by the 
eigenvalues of the matrices $A_j$ or $M_j$.}
\end{defi}

\begin{rem}\label{computationofd}
{\em 1) In the case of diagonalizable matrices $A_j$ or $M_j$ the JNF 
$J_j^n$ is completely defined by the MV $\Lambda _j^n$ and the construction 
$\Psi$ results in decreasing the biggest component of $\Lambda _j^n$ by 
$n-n_1$ to obtain $\Lambda _j^{n_1}$.  

2) For a diagonal JNF defined by a MV $\Lambda _j^n$ one has $r_j=n-m_{1,j}$ 
and $d_j=n^2-\sum _{\nu =1}^{i_j}m_{\nu ,j}^2$.

3) If $\Lambda _j^n=(n)$ and if the matrix $A_j$ or $M_j$ 
is diagonalizable, then it is scalar.}
\end{rem}

\subsection{Generic eigenvalues; non-genericity relations; the quantities 
$l$ and $\xi$\protect\label{GEPsi}}
  
We presume the necessary condition $\prod \det (C_j)=1$ (resp. 
$\sum$Tr$(c_j)=0$) to hold. This means that the eigenvalues $\sigma _{k,j}$  
(resp. $\lambda _{k,j}$) of the matrices from $C_j$ (resp. $c_j$) repeated 
with their multiplicities, satisfy the condition  
    
\begin{equation}\label{eigenv} 
\prod _{k=1}^n\prod _{j=1}^{p+1}\sigma _{k,j}=1~~{\rm resp.~~}  
\sum _{k=1}^n\sum _{j=1}^{p+1}\lambda _{k,j}=0
\end{equation} 

An equality of the form  
\[ \prod _{j=1}^{p+1}\prod _{k\in \Phi _j}\sigma _{k,j}=1~,~~{\rm resp.~~} 
\sum _{j=1}^{p+1}\sum _{k\in \Phi _j}\lambda _{k,j}=0~,\] 
is called a 
{\em non-genericity relation};  
the sets $\Phi _j$ contain one and the same number $<n$ of indices  
for all $j$. Eigenvalues satisfying none of these relations are called 
{\em generic}. Reducible  
$(p+1)$-tuples exist only for non-generic eigenvalues (a reducible 
$(p+1)$-tuple of matrices can be conjugated to a block upper-triangular form, 
its restriction to each diagonal block is such a $(p+1)$-tuple of smaller 
size, and, hence, the 
eigenvalues of each diagonal block satisfy condition (\ref{M_j}) or 
(\ref{A_j}) which is a non-genericity relation). 

\begin{rem}\label{rem1}
{\em In the case of matrices 
$A_j$, if the greatest common divisor $q$ of the multiplicities of all 
eigenvalues of all $p+1$ matrices is $>1$, then a non-genericity relation 
$(\gamma _B)$ (called the {\em basic non-genericity relation}) 
results automatically from $\sum$Tr$(c_j)=0$ when one decreases $q$ times 
the multiplicities of all eigenvalues. In  
the case of matrices $M_j$  
the equality $\prod \sigma _{k,j}=1$ implies that 
if one divides by $q$ the multiplicities of all eigenvalues, then their 
product would equal $\xi =\exp (2\pi ik/q)$, $0\leq k\leq q-1$, 
not necessarily 1. In this case a 
non-genericity relation holds exactly if $\xi$ is a non-primitive root 
of unity of order $q$. Indeed, denote by $l$ the greatest common divisor of 
$q$ and $k$. Then the product of all eigenvalues with multiplicities divided 
by $l$ equals 1 which is  the {\em basic non-genericity relation} 
$(\gamma _B)$ in the case of matrices $M_j$.}
\end{rem} 
 
\begin{defi}
{\em In the case when the basic non-genericity relation $(\gamma _B)$ 
holds eigenvalues satisfying no non-genericity relation other than 
$(\gamma _B)$ and its corollaries are called {\em relatively generic}.}
\end{defi}

The following theorem is the basic result from \cite{Ko3}, \cite{Ko4} 
and \cite{Ko5}:

\begin{tm}\label{generic}
Let $n>1$. The DSP is solvable for the conjugacy classes $C_j$ or 
$c_j$ (with generic eigenvalues,  
defining the JNFs $J_j^n$ and satisfying conditions $(\alpha _n)$ and 
$(\beta _n)$) if and only if either $\{ J_j^n\}$ satisfies 
condition $(\omega _n)$ or the construction 
$\Psi :\{ J_j^n\}\mapsto \{ J_j^{n_1}\}$ iterated as long as it is defined 
stops at a $(p+1)$-tuple $\{ J_j^{n'}\}$ either with $n'=1$ or satisfying 
condition $(\omega _{n'})$.
\end{tm}

\begin{prop}\label{psiinvar}
The construction $\Psi$ preserves the index of rigidity.
\end{prop}

The proposition is proved in \cite{Ko4}.

\begin{rem}\label{genericrem}
{\em 1) The result of the theorem does not depend on the choice one makes in 
$\Psi$ of an eigenvalue with maximal number of Jordan blocks 
(if such (a) choice(s) is (are) possible). 

2) Proposition \ref{psiinvar} implies that it suffices to check  
condition $(\alpha _{n'})$ 
for the $(p+1)$-tuple of JNFs $J_j^{n'}$ without checking  
$(\alpha _n)$ for the JNFs $J_j^n$. It does hold -- if $n'=1$, then 
$(\alpha _{n'})$ is an equality (this is the {\em rigid} case, i.e. 
$\kappa =2$). If $n'>1$ 
and condition $(\omega _{n'})$ holds for the JNFs $J_j^{n'}$, then 
$(\alpha _{n'})$ holds and is a strict inequality, see \cite{Ko3}, Theorem~9. 
Thus a posteriori one knows that it is not necessary to check condition 
$(\alpha _n)$ in Theorem~\ref{generic}.}
\end{rem}

\section{The basic result}

\subsection{The case $\kappa =0$ for diagonalizable matrices}

\begin{lm}\label{TC=Irr}
In the case $\kappa =0$ a monodromy group with a trivial centralizer and 
with relatively generic eigenvalues is irreducible.
\end{lm}

The lemma is proved in \cite{Ko5}, see part 1) of Lemma 6 there. Making use 
of the lemma we shall not distinguish solvability from weak solvability of the 
DSP in the case $\kappa =0$.

\begin{tm}\label{critnec}
In the case of matrices $M_j$, for $\kappa =0$, the conditions of 
Theorem~\ref{generic} upon 
the JNFs $J_j^n$ are 
necessary for the solvability of the DSP in the case $\kappa =0$. 
If the conjugacy 
classes $C_j$ defining the JNFs $J_j^n$ satisfy condition 
$(\beta _n)$ and do not 
satisfy condition $(\omega _n)$, then the solvability of the DSP for the 
conjugacy classes $C_j$ implies the solvability of the DSP for the 
$(p+1)$-tuple of JNFs $J_j^{n_1}=\Psi (J_j^n)$ (see Subsection~\ref{drkappa}) 
for some relatively generic eigenvalues with 
the same value of $\xi$. 
\end{tm}

The theorem is proved in Section~\ref{prcritnec}. 
In order to announce 
the basic result we need to introduce some technical notions (see 
Subsections~\ref{BTTAV} and \ref{corresp}). Therefore we 
first announce the result for the case of diagonalizable matrices which does 
not need them. 

\begin{tm}\label{diagbasicres}
1) If $\kappa =0$, if the JNFs defined by the classes $C_j$ are 
diagonal, if $q>1$, if $\xi$ is a non-primitive root of unity of 
order $q$ and if the eigenvalues of the classes $C_j$ are relatively generic, 
then the DSP is not weakly solvable for matrices $M_j$ 
(hence, not solvable either).   

2) If $\kappa =0$, if the JNFs defined by the classes $c_j$ are 
diagonal, if $q>1$ and if the eigenvalues of the classes $c_j$ are relatively 
generic, then the DSP is not weakly solvable for matrices $A_j$ 
(hence, not solvable either).
\end{tm}

A plan of the proof of the theorem is given at the end of this subsection. 

\begin{rem}\label{diagbasicresrem}
{\em It is shown in \cite{Ko5} that if the conditions of Theorem~\ref{generic} 
upon the JNFs $J_j^n$ are fulfilled and if $\xi$ is a primitive root of 
unity of order $q$, then 
the DSP is weakly solvable for matrices $M_j$ and $\kappa =0$. }
\end{rem}

In the rigid case the 
construction $\Psi$ 
stops at a $(p+1)$-tuple of one-dimensional JNFs, see Theorem~\ref{generic} 
and Remark~\ref{genericrem}, part 2). 

\begin{lm}\label{fourcases}
In the case when $\kappa =0$ and the JNFs $J_j^n$ are diagonal 
there are four possible $(p+1)$-tuples of JNFs at which $\Psi$ stops. 
Their PMVs are: 

\[ \begin{array}{llllll}{\rm Case~A)}&p=3&(d,d)&(d,d)&(d,d)&(d,d)\\
{\rm Case~B)}&p=2&(d,d,d)&(d,d,d)&(d,d,d)&\\
{\rm Case~C)}&p=2&(d,d,d,d)&(d,d,d,d)&(2d,2d)&\\
{\rm Case~D)}&p=2&(d,d,d,d,d,d)&(2d,2d,2d)&(3d,3d)&\end{array}\]
In all cases $d\in {\bf N}^*$; 
we assume that if when iterating $\Psi$ there appears a MV of the form $(n)$, 
then we delete it. 
In all four cases condition $(\omega _n)$ holds and is an equality.
\end{lm}

The lemma follows from Lemma 3 from \cite{Ko1} and from the notion of 
corresponding JNFs defined below in Subsection~\ref{corresp}.

{\bf Plan of the proof of Theorem~\ref{diagbasicres}:}
We prove part 1) first. We show that in each of the four cases A) -- D) from 
Lemma~\ref{fourcases} (and when the conditions of 1) of the theorem are 
fulfilled) the DSP is not solvable; by Lemma~\ref{TC=Irr} it is not weakly 
solvable either. This is done in 
Sections~\ref{seccaseA}, \ref{seccaseB}, \ref{seccaseC} and \ref{seccaseD}, 
one case per section. Section~\ref{seccaseA} is the 
longest and the most important of them because in the other three cases the 
proof is reduced to the one in Case A). 

Theorem~\ref{critnec} and  
Lemma~\ref{fourcases} imply that in all 
possible cases covered by Theorem~\ref{diagbasicres} the DSP is not weakly 
solvable.
Part 2) of the theorem is proved in Section~\ref{matricesA_j} using part 1).

\subsection{The basic technical tool\protect\label{BTTAV}}

\begin{defi}
{\em Call {\em basic technical tool} the way described 
below to deform analytically a $(p+1)$-tuple of matrices $A_j$ satisfying 
(\ref{A_j}) or of matrices $M_j$ satisfying (\ref{M_j}) with a 
{\em trivial centralizer}.}
\end{defi} 

In the case of matrices $A_j$ set 
$A_j=Q_j^{-1}G_jQ_j$, 
$G_j$ being Jordan matrices. Look for matrices $\tilde{A}_j$ of the form 
$\tilde{A}_j=(I+\sum _{i=1}^s\varepsilon _iX_{j,i}(\varepsilon ))^{-1}
Q_j^{-1}
(G_j+\sum _{i=1}^s\varepsilon _iV_{j,i}(\varepsilon ))Q_j
(I+\sum _{i=1}^s\varepsilon _iX_{j,i}(\varepsilon ))$  
where 
$\varepsilon =(\varepsilon _1,\ldots ,\varepsilon _s)\in ({\bf C}^s,0)$ and 
$V_{j,i}(\varepsilon )$ are given 
matrices analytic 
in $\varepsilon$. One chooses $V_{j,i}$ such that  
tr$(\sum _{j=1}^{p+1}\sum _{i=1}^s\varepsilon _iV_{j,i}(\varepsilon ))
\equiv 0$ identically in $\varepsilon$. One often has $s=1$ and 
$V_{j,1}$ are such that the 
eigenvalues of the $(p+1)$-tuple of matrices $\tilde{A}_j$ are generic for 
$\varepsilon \neq 0$. Often one has $V_{j,i}\equiv 0$ for all indices $j$ 
but one, i.e. all matrices $A_j$ but one remain within their conjugacy 
classes.

In the case of $(p+1)$-tuples of matrices $M_j^1$ with a 
trivial centralizer look for 
$M_j$ of the form 

\begin{equation}\label{BTTM} 
M_j=(I+\sum _{i=1}^s\varepsilon _iX_{j,i}(\varepsilon ))^{-1}(M_j^1+
\sum _{i=1}^s\varepsilon _iN_{j,i}(\varepsilon ))(I+
\sum _{i=1}^s\varepsilon _iX_{j,i}(\varepsilon ))
\end{equation} 
where the given matrices $N_{j,i}$ are analytic in 
$\varepsilon \in ({\bf C}^s,0)$ and one looks for $X_{j,i}$ analytic in 
$\varepsilon$. Like in the case of matrices $A_j$ one can set 
$M_j^1=Q_j^{-1}G_jQ_j$, $N_{j,i}=Q_j^{-1}V_{j,i}Q_j$. For both cases the 
existence of the matrices $X_{j,i}$ analytic in $\varepsilon$  
is proved in \cite{Ko4}.

\subsection{Correspondence between Jordan normal forms\protect\label{corresp}}

\begin{defi}
{\em For a given JNF $J^n=\{ b_{i,k}\}$ define its {\em corresponding} 
diagonal JNF ${J'}^n$. A diagonal JNF is  
a partition of $n$ defined by the multiplicities of the eigenvalues. 
For each $k$ fixed the collection $\{ b_{i,k}\}$ is a partition ${\cal P}_k$ 
of $\sum _{i\in I_k}b_{i,k}$. The diagonal JNF  
${J'}^n$ is the disjoint sum of the partitions dual to ${\cal P}_k$.} 
\end{defi}

\begin{ex}
{\em Consider the JNF $J^{17}=\{ \{ 6,4,3\} \{ 3,1\} \}$, i.e. 
with two eigenvalues, 
the first with three Jordan blocks of sizes 
6,4,3  and the second with two blocks of sizes 3,1. The partition of 13 
dual to (6,4,3) is (3,3,3,2,1,1), the one of 4 dual to  
(3,1) is (2,1,1). Hence, the diagonal JNF 
corresponding to $J^{17}$ is defined by the MV (3,3,3,2,2,1,1,1,1) (in  
decreasing order of the multiplicities).}
\end{ex}

\begin{prop}\label{psicorr}
Consider a JNF $J^n$ and its corresponding diagonal JNF ${J'}^n$ defined 
by a MV $\Lambda =(m_1,\ldots ,m_{\nu})$, 
$m_1\geq \ldots \geq m_{\nu}$. Choose an 
eigenvalue of $J^n$ with maximal number $n-r(J^n)$ of Jordan blocks and 
decrease the sizes of the $k'$ smallest of these blocks by 1, 
$k'\leq n-r(J^n)$ -- this defines a new JNF $J^{n-k'}$. 
Set $\Lambda _*=(m_1-k',m_2,\ldots ,m_{\nu})$. Then the MV $\Lambda _*$ 
defines a diagonal JNF corresponding to $J^{n-k'}$.
\end{prop}

\begin{cor}\label{corpsicorr}
The $(p+1)$-tuples of JNFs $J_j^n$ ${J_j'}^n$ where for each $j$ 
$J_j^n$ corresponds to ${J_j'}^n$ satisfy or not the conditions of 
Theorem~\ref{generic} simultaneously.
\end{cor} 

The propositions and corollary from this subsection are proved in \cite{Ko4}.

\begin{prop}\label{rd}  
1) If the JNF ${J'}^n$ corresponds to the JNF $J^n$, then 
$r(J^n)=r({J'}^n)$ and $d(J^n)=d({J'}^n)$.

2) To each diagonal JNF there corresponds a unique JNF with a single 
eigenvalue.
\end{prop}

\begin{rem}\label{deformation}
{\em Denote by $G$ a Jordan matrix and by $G'$ a diagonal matrix defined as 
follows: the diagonal entries of $G'$ in the last but $s$ positions of the 
Jordan 
blocks of $G$ with given eigenvalue $\lambda$ are equal among themselves 
and different from the ones in the 
last but $m$ positions for $m\neq s$, $m,s\in {\bf N}^*$. Then the matrix 
$G+\varepsilon G'$, $0\neq \varepsilon \in ({\bf C},0)$ is diagonalizable 
and its JNF is the diagonal JNF corresponding to $J(G)$ (the poof can be found 
in \cite{Ko4}). Hence, if one applies the basic technical tool with 
$s=1$ and $G_j$, $V_{j,1}$ playing the roles respectively of $G$, $G'$, then 
one sees that the weak solvability of the DSP for 
matrices $A_j$ or $M_j$ with given JNFs $J_j^n$ implies the one for 
diagonal JNFs corresponding to $J_j^n$ and for nearby eigenvalues.}
\end{rem}

\subsection{The result in the general case\protect\label{resgencase}}

\begin{defi}
{\em We say that the conjugacy class 
$C$ is {\em continuously deformed} into the class $C'$ if either the 
classes $C$, $C'$ are like the ones of the matrices $G$, 
$G+\varepsilon G'$ from Remark~\ref{deformation} or $C'$ is just another 
conjugacy class defining the same JNF as $C$. We say that the $(p+1)$-tuple 
of conjugacy classes $C_j$ is continuously deformed into the $(p+1)$-tuple of 
conjugacy classes $C_j'$ if each class $C_j$ is continuously deformed into the 
corresponding class $C_j'$ and the eigenvalues of the first $(p+1)$-tuple 
are homotopic to the ones of the second $(p+1)$-tuple. Throughout the homotopy 
there holds condition (\ref{eigenv}) and the MVs remain the same.}
\end{defi}

\begin{ex}\label{twodifferentq}
{\em Consider the triple of conjugacy classes $C_1$, $C_2$, $C_3$ of size 12 
each with a single eigenvalue $\lambda _j$ 
and with Jordan blocks of equal size $l_j$: 
$(\lambda _1,\lambda _2,\lambda _3)=(i,1,1)$, $(l_1,l_2,l_3)=(2,3,6)$.
For these eigenvalues one has $q=12$, $\xi =i$ which is not a primitive 
root of unity of order 12. One has $l=3$. 
The basic non-genericity relation $(\gamma _B)$ 
is obtained by dividing the multiplicities of all eigenvalues by 3. The 
eigenvalues are relatively generic.

To the triple of JNFs defined by the  
conjugacy classes $C_j$ there corresponds the triple of diagonal JNFs 
defined by the PMV $(6,6)$, $(4,4,4)$, $(2,2,2,2,2,2)$. For this PMV 
one has $q=2$ and by continuous deformation of the conjugacy classes $C_j$ 
into diagonal ones with the above PMV one obtains $\xi =-1$ which is a 
primitive root of unity of order 2. (Indeed, for the classes $C_j$ the product 
of the eigenvalues repeated each with the half of its multiplicity 
equals $-1$ which remains unchanged throughout the continuous deformation.)}
\end{ex} 

\begin{defi}
{\em Denote by $d$  
the greatest common divisor of all quantities 
$\Sigma _{j,m}(\sigma )$ where $\Sigma _{j,m}(\sigma )$ is the number of 
Jordan blocks of size $m$ of a given matrix $M_j$ or $A_j$ and with eigenvalue 
$\sigma$. It is true that $d$ divides $q$ and that $q$ divides $n$.}
\end{defi}

\begin{rem}
{\em The quantity $q$ does not increase under continuous deformations 
like in the above example. If one deforms continuously 
the conjugacy classes so that the eigenvalues of $C'$ 
be ``as generic as possible'' (i.e. satisfying only these non-genericity 
relations which are not destroyed by continuous deformations like 
the above ones), then one has $q=d$.}
\end{rem}

\begin{tm}\label{maintm}
Suppose that 

1) the conjugacy classes of the matrices $A_j$ or $M_j$ verify the 
conditions of Theorem~\ref{generic};

2) they are continuously deformed into a $(p+1)$-tuple of conjugacy classes 
defining diagonal JNFs with $q=d>1$, with relatively generic eigenvalues 
and in the case of matrices 
$M_j$ with $\xi$ being a non-primitive root of unity of order $q$; 

3) one has $\kappa =0$.  

Then for such conjugacy classes the DSP is not weakly solvable.
\end{tm}

{\bf Proof:}
Suppose that there exists a $(p+1)$-tuple 
of matrices $M_j$ with trivial centarlizer which satisfies conditions 1), 2) 
and 3). Applying the basic technical tool with $l=1$ and $G_j$, $V_{j,1}$ like 
in Remark~\ref{deformation}, one obtains the existence of a $(p+1)$-tuple of 
diagonalizable matrices $M_j$ with a trivial centralizer, with relatively 
generic eigenvalues, with $\kappa =0$ 
and with $\xi$ being a non-primitive root of unity of order $q$ which 
contradicts Theorem~\ref{diagbasicres}.\hspace{1cm}$\Box$

\section{Proof of Theorem \protect\ref{critnec}\protect\label{prcritnec}}
\subsection{The proof itself\protect\label{pritself}}

\begin{defi}
{\em A regular singular point of a linear system of ordinary differential 
equations is called {\em apparent} if 
its local monodromy is trivial.}
\end{defi}

\begin{lm}\label{addapp}
Any monodromy group can be realized by a Fuchsian system on ${\bf C}P^1$ with 
at most one additional apparent singularity at a point $a_{p+2}$ 
which can be chosen 
arbitrarily; for the eigenvalues $\lambda _{k,j}$ of 
the matrices-residua $A_j$, 
$j=1,\ldots, p+1$ one has {\rm Re}$\lambda _{k,j}\in [0,1)$; 
one has $J(A_j)=J(M_j)$ for $j=1,\ldots ,p+1$, $M_j$ being the monodromy 
operators. 
\end{lm}

The lemmas from this subsection except Lemmas~\ref{norm} and 
\ref{subreprtrivcentr} are proved in the subsequent ones (one proof 
per subsection). In what follows the points $a_1$, 
$\ldots$, $a_{p+2}$ are fixed. 

\begin{defi}
{\em A Fuchsian system belongs to the {\em class N} if it has poles at the 
points $a_j$ the one at 
$a_{p+2}$ being an apparent singularity, if its monodromy group is 
irreducible, and if at $a_{p+2}$ the Laurent series 
expansion of the system looks like this:}

\begin{equation}\label{classN} 
\dot{X}=(A_{p+2}/(t-a_{p+2})+B(t-a_{p+2}))X
\end{equation}
{\em where $A_{p+2}=$diag$(\mu _1,\ldots ,\mu _n)$, $\mu _j\in {\bf Z}$, 
$\mu _1\geq \ldots \geq \mu _n$. 

Denote by ord$u$ the order of the zero at $a_{p+2}$ of the germ of 
holomorphic function $u$. A class N Fuchsian system is called {\em normalized} 
if for $i<j$ one has ord$B_{i,j}\geq \mu _i-\mu _j$.}
\end{defi}

\begin{lm}\label{norm}
If one has  
$A_{p+2}=${\rm diag}$(\mu _1,\ldots ,\mu _n)$, $\mu _j\in {\bf Z}$, 
$\mu _1\geq \ldots \geq \mu _n$, and if one has for $i<j$ 
{\rm ord}$B_{i,j}\geq \mu _i-\mu _j$ for $B$ defined by (\ref{classN}), 
then the singularity at $a_{p+2}$ is apparent.
\end{lm}

Indeed, the following change of variables brings the 
system locally, at $a_{p+2}$, to a system without a pole at $a_{p+2}$ (hence, 
the local monodromy at $a_{p+2}$ is trivial):

\begin{equation}\label{shift}
X\mapsto (t-a_{p+2})^{{\rm diag}(\mu _1,\ldots ,\mu _n)}X
\end{equation}
 
\begin{lm}\label{leqp}
For a normalized class N Fuchsian system one has $\mu _i-\mu _{i+1}\leq p$ 
for $i=1,\ldots ,n-1$.
\end{lm}

\begin{defi}
{\em Set $\sigma =(\mu _1+\ldots +\mu _n)/n$ ({\em mean value}) and 
$\delta =((\mu _1-\sigma )^2+\ldots +(\mu _n-\sigma )^2)/n$ ({\em dispersion} 
of the numbers $\mu _i$).}
\end{defi}

\begin{lm}\label{decrease}
The monodromy group of a 
non-normalized class N Fuchsian system can be realized 
by a normalized class N Fuchsian system with the same conjugacy classes of 
the matrices $A_1$, $\ldots$, $A_{p+1}$, with the same mean value and with a 
smaller dispersion of the numbers $\mu _i$.
\end{lm}
 
Suppose that for $\kappa =0$ and for 
given diagonal conjugacy classes with relatively generic 
eigenvalues and not satisfying condition $(\omega _n)$ 
there exists a monodromy group with a trivial centralizer (hence, 
irreducible by Lemma~\ref{TC=Irr}). Then for almost all relatively generic 
eigenvalues with the same value of $\xi$ there exist irreducible monodromy 
groups with such JNFs. Indeed, applying the basic technical tool, one can 
deform the given monodromy group into one with any nearby relatively generic 
eigenvalues and the same JNFs of the matrices $M_j$. Moreover, the 
deformation can be chosen such that the new matrices $M_j$ will be 
diagonalizable and defining the JNFs corresponding to the initial ones. 

The set ${\cal M}$ of such 
monodromy groups is constructible and such is its projection ${\cal V}$ 
on the set of 
eigenvalues ${\cal W}$, i.e. ${\cal V}$ is an everywhere dense 
constructible subset of ${\cal W}$.

Lemmas \ref{addapp}, \ref{leqp} and \ref{decrease} imply that for 
given conjugacy classes $C_j$ of $M_1$, $\ldots$, $M_{p+1}$ there exist 
finitely many sets $\Gamma _i$ 
of eigenvalues $\mu _k=\lambda _{k,p+2}$ such that the monodromy 
group can be realized by a normalized class N Fuchsian system with such 
eigenvalues of $A_{p+2}$; for $j\leq p+1$ the eigenvalues 
$\lambda _{k,j}$ are uniquely defined by the classes $C_j$, see 
Lemma~\ref{addapp}. 

Consider $gl(n,{\bf C})^{p+1}$ as the space of $(p+2)$-tuples of matrices 
$A_j$ whose sum is 0. Denote by ${\cal G}_i$ its subsets such that $A_{p+2}$ 
is diagonal, with eigenvalues $\mu _k\in \Gamma _i$, and for $i<j$ 
there holds the condition ord$B_{i,j}\geq \mu _i-\mu _j$ for $B$ defined by 
(\ref{classN}) (recall that the poles $a_j$ are fixed). Hence, the sets 
${\cal G}_i$ are constructible. 

A point from ${\cal G}_i$ defines a Fuchsian system (S). Fix a base point $a$ 
different from the points $a_j$ and define the monodromy operators of the 
system with initial data $X|_{t=a}=I$. The map which maps the 
matrices-residua $A_1$, $\ldots$, $A_{p+2}$ into 
the $(p+1)$-tuple of monodromy operators of system (S) is a map  
$\chi _i :{\cal G}_i\rightarrow {\cal M}$. 

For each point from ${\cal M}$ 
there exists at least one $i$ such that the point has a preimage in 
${\cal G}_i$ under $\chi _i$. 
This means that there exists a point from ${\cal M}$ such that some 
neighbourhood of his is covered by $\chi _i({\cal G}_i)$ for some $i$; 
we set $i=1$. Indeed, the constructible set 
${\cal M}$ cannot be locally 
covered by a finite number of analytic sets of lower dimension. This and 
the irreducibility of ${\cal M}$ implies 
that the set $\chi _1({\cal G}_1)$ is dense in ${\cal M}$. 

\begin{lm}\label{fivecond}
Suppose that 

A) the matrices-residua $A_1$, $\ldots$, $A_{p+1}$ of a normalized class N 
Fuchsian system are diagonalizable, with generic eigenvalues; 

B) their $(p+1)$-tuple is irreducible; 

C) none of these matrices has eigenvalues differing by a non-zero integer and  
each of them has a single integer eigenvalue $\lambda _j$ whose multiplicity 
is a (the) greatest one (hence, each 
monodromy operator $M_j$ has an eigenvalue $\sigma _j=1$);

D) all non-genericity relations satisfied by the 
eigenvalues of the monodromy operators $M_j$ result from two relations, the 
first of which is the basic one $(\gamma _B)$ the second being 

\[ \sigma _1\ldots \sigma _{p+1}=1~~~~~~~~~~~~~~(\gamma _0)\]

E) one has $\lambda _j>0$ and 
$\lambda _1+\ldots +\lambda _{p+1}>(n^2+n)\mu$ with 
$\mu =\max (|\mu _1|,|\mu _n|)$. 

F) the monodromy group can be analytically deformed into an 
irreducible one for nearby relatively  
generic eigenvalues and with the same JNFs of the matrices 
$M_j$.

G) Condition $(\omega _n)$ does not hold for the matrices $M_j$. 

Then the monodromy group of the Fuchsian system is 
with trivial centralizer.
\end{lm}

The projection ${\cal P}_1$ of the set ${\cal G}_1$ on the space ${\bf C}^s$ 
of eigenvalues $\lambda _{k,j}$ ($s$ depends on their multiplicities)  
is a constructible set. If ${\cal P}_1$  
does not contain a point satisfying conditions C), D) and E) of the lemma, 
then codim$_{{\bf C}^s}{\cal P}_1>0$, 
hence, $\chi _1({\cal G}_1)$ cannot be dense in ${\cal M}$.

\begin{lm}\label{theMGisreducible}
The monodromy group of system (\ref{Fuchs}) with eigenvalues defined as 
in Lemma~\ref{fivecond} can be conjugated to the form 
$\left( \begin{array}{cc}\Phi &\ast \\
0&I\end{array}\right)$ where $\Phi$ is $n_1\times n_1$.
\end{lm}

The subrepresentation $\Phi$ can be reducible. The following 
lemma is proved in \cite{Ko4}.

\begin{lm}\label{subreprtrivcentr}
The centralizer ${\cal Z}(\Phi )$ of the subrepresentation $\Phi$ is trivial.
\end{lm}

Thus the existence of an irreducible representation of rank $n$ for which 
condition $(\omega _n)$ does not hold implies the existence of the 
representation $\Phi$ of rank $n_1$ and with trivial centralizer. The JNFs 
defined by the matrices from $\Phi$ are obtained from the initial $(p+1)$ 
JNFs by applying the map $\Psi$. One can 
deform the eigenvalues of $\Phi$ so that they become relatively generic. For 
such eigenvalues the deformed representation $\Phi$ is irreducible, see 
Lemma~\ref{TC=Irr}. If 
$\Phi$ satisfies condition $(\omega _{n_1})$, then we are done. 
If not, then we 
continue iterating $\Psi$. In the end we stop at a representation of rank 
$n'$ satisfying 
condition $(\omega _{n'})$. It is impossible to obtain a representation of 
rank 1 because its index of rigidity is 2, see Proposition~\ref{psiinvar}.

The eigenvalues of the representation $\Phi$ define the same value of $\xi$ 
as the ones of the initial representation. Indeed, the eigenvalues from the 
initial one which are not in $\Phi$ equal 1.~~~~$\Box$

\subsection{Proof of Lemma~\protect\ref{addapp}}

It is shown in \cite{P} that any monodromy group can be realized by a 
regular system on ${\bf C}P^1$ which is Fuchsian at all poles but one. So one 
can add a $(p+2)$-nd monodromy operator equal to $I$ to the initial operators 
$M_j$ assuming that the system realizing this monodromy group has not $p+1$ 
but $p+2$ poles. Applying the result from \cite{P} (reproved in \cite{ArIl}, 
p. 131)  
one obtains a regular 
system (S) with the given monodromy group which is Fuchsian at 
$a_1$, $\ldots$, 
$a_{p+1}$ and which has a regular apparent singularity at $a_{p+2}$. The 
point $a_{p+2}\neq a_j$, $j\leq p+1$, is chosen arbitrarily and  
the JNFs of the matrices 
$A_j$ are the same as the ones of the corresponding monodromy operators 
$M_j$ for $j=1,\ldots ,p+1$. Moreover, Re$\lambda _{k,j}\in [0,1)$. 

\begin{rem}
{\em In \cite{P} an attempt is made to prove that every monodromy group can be 
realized by a Fuchsian system on ${\bf C}P^1$ (without apparent 
singularities). This is one of the versions of the Riemann-Hilbert problem 
and the answer to it is negative, see \cite{Bo1}. We are referring above 
to the correct part of the attempt from \cite{P} to prove the Riemann-Hilbert 
problem. See \cite{ArIl} pp. 130 -- 135 as well.}
\end{rem}

\noindent Make the singularity at $a_{p+2}$  
Fuchsian. Fix a matrix solution to system (\ref{Fuchs}) with 
$\det X\not\equiv 0$. Its regularity and the triviality of the monodromy at 
$a_{p+2}$ imply that it is meromorphic at $a_{p+2}$. 

\begin{lm}\label{Souvage}
(A. Souvage) A meromorphic mapping from ${\bf C}^n$ to ${\bf C}^n$ with a pole 
at $a_{p+2}$ and nondegenerate for $t\neq a_{p+2}$ 
can be represented in the form $PH(t-a_{p+2})^D$ where $D$ is a 
diagonal matrix with integer entries, $H$ is holomorphic and holomorphically 
invertible at $a_{p+2}$ and the entries of the matrix $P$ are polynomials 
in $1/(t-a_{p+2})$, $\det P\equiv$const$\neq 0$. 
\end{lm} 

Perform in system (S) the change $X\mapsto P^{-1}X$. This change leaves the 
system Fuchsian at $a_1$, $\ldots$, $a_{p+1}$ and regular at $a_{p+2}$ without
introducing new singular points. At $a_{p+2}$ the new system is Fuchsian. 
Indeed, the matrix $(t-a_{p+2})^D$ is a solution to the system (Fuchsian 
at $a_{p+2}$) $\dot{X}=(D/(t-a_{p+2}))X$. The change of variables 
$X\mapsto HX$ leaves the 
latter system Fuchsian at $a_{p+2}$ (the system becomes 
$\dot{X}=(-H^{-1}\dot{H}+H^{-1}(D/(t-a_{p+2}))H)X$).~~~~$\Box$

\subsection{Proof of Lemma \protect\ref{leqp}}

$1^0$. The matrix $B$ defined by equation (\ref{classN}) admits the 
Taylor series expansion 
$B=B_0+(t-a_{p+2})B_1+(t-a_{p+2})^2B_2+\ldots$. A 
direct computation shows that 
$B_{\nu}=-\sum _{j=1}^{p+1}A_j/(a_j-a_{p+2})^{\nu}$.
Suppose that for some $i_0$ ($1\leq i_0\leq n-1$) one has 
$\mu _{i_0}-\mu _{i_0+1}\geq p+1$. Then for $i\leq i_0$, $k\geq i_0+1$ one has 
$\mu _i-\mu _k\geq p+1$. 

$2^0$. Hence, all matrix entries $A_{j;i,k}$ with $j\leq p+1$ and $i,k$ like 
in $1^0$ must be 0. Indeed, for each such $i,k$ fixed 
the system of linear equations 
$B_{\nu ;i,k}=0$, $\nu =1,\ldots ,p+1$  
with unknown variables the entries $A_{j;i,k}$ implies $A_{j;i,k}=0$ because 
it is of rank $p+1$ (its determinant is the Vandermonde one 
$W(1/(a_1-a_{p+2}),\ldots ,1/(a_{p+1}-a_{p+2}))$ and for $j_1\neq j_2$ one 
has $a_{j_1}\neq a_{j_2}$). 

This means that the matrices-residua $A_1$, $\ldots$, $A_{p+1}$ 
are block lower-triangular, with diagonal 
blocks of sizes $i_0$ and $n-i_0$. Hence, so are the monodromy operators, i.e. 
the monodromy group is reducible and the system is not from the 
class N.~~~~$\Box$

\subsection{Proof of Lemma \protect\ref{decrease}}

$1^0$. Recall that the matrix $B$ was defined by equation (\ref{classN}). 
Assume for simplicity that $a_{p+2}=0$. 
For $i<j$ find an entry $B_{i,j}$ with smallest value of 
$m:=-$ord$B_{i,j}-\mu _j+\mu _i$. Hence, $m>0$. If there are several 
possible choices, then we choose among them one with minimal value of 
$j-i$. Set $B_{i,j}=bt^g+
o(|t|^g)$, $b\neq 0$ (hence, $g=$ord$B_{i,j}$). 

$2^0$. Consider the change of variables $X\mapsto WX$ with 
$W=I+(\mu _j-\mu _i+g)E_{j,i}/bt^m$. It is holomorphic for $t\neq 0$, with 
$\det W\equiv 1$, hence, it preserves the conjugacy classes of the 
residua $A_1$, $\ldots$, $A_{p+1}$ the system remaining Fuchsian there. 
At $a_{p+2}$ the new residuum is 
lower-triangular, with diagonal entries equal to 
$\mu _1,\ldots ,\mu _{i-1},\mu _j+g,\mu _{i+1},\ldots ,\mu _{j-1},\mu _i-g,
\mu _{j+1},\ldots ,\mu _n$. 
The singularity at $a_{p+2}$, in general, is no longer Fuchsian, but the 
order of the pole at $a_{p+2}$ is $\leq m$; equality is possible only in 
position $(j,i)$. This follows from rule (\ref{gauge}) (the reader is invited 
to check the claim). 

Except on the diagonal poles of order $>1$ at 0 can appear only in the entries 
$(j,1)$, $(j,2)$, $\ldots$, $(j,i)$, $(j+1,i)$, $(j+2,i)$, $\ldots$, $(n,i)$, 
see the choice of $B_{i,j}$ in $1^0$. 

$3^0$. One deletes the polar terms below the diagonal by a change 
$X\mapsto VX$, $V=I+V'$ where each entry $V'_{k,\nu }$ 
of $V'$ is a suitably chosen polynomial $p_{k,\nu }$ of 
$1/t$, the non-zero entries being in the positions cited at the end of $2^0$. 
The degree of the polynomial $p_{k,\nu }$ is equal to the order of the pole 
in position $(k,\nu)$ which has to disappear. We leave for the reader the 
proof that such a choice of the polynomials $p_{k,\nu }$ is really possible.

$4^0$. As a result of the changes from $2^0$ and $3^0$ the system remains 
Fuchsian at $a_j$ for $j\leq p+1$ and the conjugacy classes of its residua 
do not change because the matrix $V$ is  
holomorphic for $t\neq 0$ and $\det V\equiv 1$. The system 
remains Fuchsian at 0 as well and the eigenvalues of $A_{p+2}$ change as 
follows: $\mu _i\mapsto \mu _i-g$, $\mu _j\mapsto \mu _j+g$, the rest of the 
eigenvalues remain the same. (One should rearrange after this the 
eigenvalues $\mu _i$ in 
decreasing order by conjugating with a constant permutation matrix.) 
One checks directly that as a result of the change of the eigenvalues 
$\mu _i$ the mean value $\sigma$ remains the same  
whereas $\delta$ decreases. ~~~~$\Box$

\subsection{Proof of Lemma \protect\ref{fivecond}}

$1^0$. Suppose that the centralizer ${\cal Z}$ is nontrivial. 
Hence, it contains either a diagonalizable matrix $D$ with exactly 
two different eigenvalues or 
a nilpotent matrix $N\neq 0$ such that $N^2=0$.

$2^0$. Suppose that 
$D=\left( \begin{array}{cc}\alpha I&0\\0&\beta I\end{array}\right)
\in {\cal Z}$ with 
diagonal blocks of sizes $l'$ and $n-l'$ and with $\alpha \neq \beta$. 
Then the matrices $M_j$ are block-diagonal with the same sizes of the 
diagonal blocks and the monodromy group is a 
direct sum. This follows from $[M_j,D]=0$. Denote 
the two diagonal blocks of $M_j$ by $S_j$ and $T_j$ ($S_j$ is $l'\times l'$).

Hence, there are two subspaces of the solution space (${\cal X}_1$ and 
${\cal X}_2$) which are invariant for the monodromy group and whose direct sum 
is the solution space. Denote by $C_j'$, $C_j''$ the conjugacy classes of the 
matrices $S_j$ and $T_j$.

$3^0$. Use a result from \cite{Bo1} (see Lemma 3.6 there): 

\begin{lm}\label{Z}
The sum of the eigenvalues $\lambda _{k,j}$ of the matrices-residua $A_j$ 
corresponding to an invariant subspace of the monodromy group is a 
non-positive integer.
\end{lm}

\begin{rem}\label{warning}
{\em 1) Condition C) and Remark~\ref{exp} imply that the 
equality $\exp (2\pi i\lambda _{k,j})=\sigma _{k,j}$ defines 
(for $j\leq p+1$ fixed) a bijection 
between the eigenvalues $\sigma _{k,j}$ and the eigenvalues $\lambda _{k,j}$ 
modulo permutation of equal eigenvalues. For $j=p+2$ this is false (recall 
that $\lambda _{k,p+2}=\mu _k\in {\bf Z}$, 
$\sigma _{1,p+2}=\ldots =\sigma _{n,p+2}=1$). 

2) When defining the sets of eigenvalues $\lambda _{k,j}$ corresponding to the 
subspaces ${\cal X}_1$ and ${\cal X}_2$ it is true only for $j\leq p+1$ 
but not for $j=p+2$ that these sets are complementary to one another, i.e. 
one and the same eigenvalue $\lambda _{k,p+2}=\mu _k$ might appear in both 
sums while another one might appear in none of them. 

Indeed, present  
the eigenvalues $\lambda _{k,j}$ in the form 
$\varphi _{k,j}+\rho _{k,j}$ with $\varphi _{k,j}\in {\bf Z}$, 
Re$\rho _{k,j}\in [0,1)$ (this presentation is unique). The numbers 
$\varphi _{k,j}$ have the meaning of valuations on the solution 
subspace on which 
the monodromy operator $M_j$ acts with a single 
eigenvalue $\exp (2\pi i\rho _{k,j})$, 
see the details in \cite{Bo1} (Definition 2.3 etc.). 
 
At $a_{p+2}$ one has $\rho _{k,p+2}=0$, $\varphi _{k,p+2}=\mu _k$. 
Thus if a vector-column solution $\tilde{X}'\in {\cal X}_1$ of 
system (\ref{Fuchs}) has an 
expansion at $a_{p+2}$ into a Laurent series  
$v_1(t-a_{p+2})^{\mu _{i_1}}+v_2(t-a_{p+2})^{\mu _{i_2}}+
o((t-a_{p+2})^{\mu _{i_2}})$, with $\mu _{i_1}<\mu _{i_2}$ and 
$0\neq v_i\in {\bf C}^n$, then it is 
$\mu _{i_1}$ that participates in the sum of eigenvalues $\lambda _{k,j}$ 
corresponding to 
${\cal X}_1$ because this is the valuation of $\tilde{X}'$ at $a_{p+2}$. 

If a solution $\tilde{X}''\in {\cal X}_2$ equals 
$cv_1(t-a_{p+2})^{\mu _{i_1}}+dv_2(t-a_{p+2})^{\mu _{i_2}}+
o((t-a_{p+2})^{\mu _{i_2}})$, $c,d\in {\bf C}^*$, $c\neq d$, then it is again 
$\mu _{i_1}$ that participates in the sum corresponding to 
${\cal X}_2$. The number $\mu _{i_2}$ is a valuation of the 
solution $c\tilde{X}'-\tilde{X}''$ which might be neither in ${\cal X}_1$ nor 
in ${\cal X}_2$, therefore $\mu _{i_2}$ might appear in neither of the 
two sums.
For $j\leq p+1$ there is no such ambiguity due to condition C), i.e. to 
each eigenvalue of the monodromy operator $M_j$ there corresponds a single 
valuation on the corresponding solution subspace.}
\end{rem}

$4^0$. {\em Lemma~\ref{Z} and conditions D) and E) imply that if the 
monodromy group is a 
direct sum, then equal eigenvalues of the matrices 
$S_j$ and $T_j$ have proportional multiplicities.} 

Indeed, denote by $\Xi$, $\Theta$ the sets 
of eigenvalues $\sigma _{k,j}$, $j\leq p+1$ participating respectively 
in $(\gamma _B)$, $(\gamma _0)$ and by $\Xi '$, $\Theta '$ the sums of their 
respective eigenvalues $\lambda _{k,j}$. 
Hence, the sums of eigenvalues of 
the matrices $A_1$, $\ldots$, $A_{p+2}$ relative to the solution 
subspaces ${\cal X}_1$ and 
${\cal X}_2$ are both of the form $\phi _i:=a_i\Xi '+b_i\Theta '+\Delta _i$, 
$a_i\in {\bf N}$, $b_i\in {\bf Z}$, $b_1+b_2=0$ where $\Delta _1$ 
(resp. $\Delta _2$) is the sum of some $l'$ (resp. $n-l'$) eigenvalues 
$\lambda _{k,p+2}=\mu _k$ (see Remark~\ref{warning}); hence, 
$|\Delta _i|\leq n\mu$. 

One has $a_i\leq n$ (evident), and $|\Xi '|<n\mu$ (because the sum of all 
eigenvalues $\lambda _{k,j}$ (which is 0) is of the form 
$g\Xi '+\sum _{k=1}^n\mu _k$ with $g\in {\bf N}$, $1<g<n$; hence, 
$|\Xi '|\leq n\mu /g<n\mu$). 

If $b_1>0$, then 
$\phi _1\geq b_1(n^2+n)\mu -a_1|\Xi '|-|\Delta _1|>(n^2+n)\mu -n^2\mu -
n\mu >0$.  
This contradicts Lemma~\ref{Z}. Hence, $b_1\leq 0$. In the same way 
$b_2\leq 0$. Hence, $b_1=b_2=0$. This means that equal eigenvalues of 
the blocks $S_j$ and $T_j$ have proportional multiplicities.

$5^0$. The monodromy group of a Fuchsian system satisfying the condition 
$b_1=b_2=0$, see $4^0$, cannot 
be analytically deformed into an irreducible one for nearby relatively  
generic eigenvalues and with the same Jordan normal forms of the matrices 
$M_j$; this contradicts condition F). 

Indeed, suppose that there exists such a 
deformation analytic in $\varepsilon \in ({\bf C},0)$ (i.e. for 
almost all values of $\varepsilon \neq 0$ the $(p+1)$-tuple is irreducible). 
For the $(p+1)$-tuple before the deformation the 
multiplicities of the equal eigenvalues $\sigma _{k,j}$ of the two diagonal 
blocks $S_j$ and $T_j$ are proportional for all $j$. This means that for all 
$j$ one has $d(C_j')=({l'}^2/n^2)d(C_j)$, $d(C_j'')=((n-l')^2/n^2)d(C_j)$. 
Indeed, if a diagonal JNF is defined by the PMV $(m_1,\ldots ,m_s)$, 
then a conjugacy class defining such a JNF is of dimension  
$n^2-\sum _{i=1}^s(m_i)^2$. 
Hence, $d(C_1'')+\ldots +d(C_{p+1}'')=2(n-l')^2$, 
$d(C_1')+\ldots +d(C_{p+1}')=2{l'}^2$ (this follows from the proportional 
multiplicities) 
and for the representations ${\cal M}'$, ${\cal M}''$ defined by the 
matrices $S_j$, $T_j$ one has 

\begin{equation}\label{Ext1}
{\rm Ext}^1({\cal M}',{\cal M}'')={\rm Ext}^1({\cal M}'',{\cal M}')=0
\end{equation}

$6^0$. When one deforms analytically a $(p+1)$-tuple into a nearby one 
(see the 
basic technical tool) one can express the deformation as a superposition of 
two deformations -- of a change of the eigenvalues (see the matrices 
$N_{j,i}(\varepsilon )$ in (\ref{BTTM})) and of a conjugation (see the 
matrices $X_{j,i}(\varepsilon )$ there). One can choose the matrices $N_{j,i}$ 
to be polynomials of the matrices $M_j$, i.e. block-diagonal, with diagonal 
blocks of sizes $l'$ and $n-l'$. Hence, the two non-diagonal blocks of the 
matrices change (in first approximation w.r.t. $\varepsilon$) only 
as a result of the conjugation. 

Condition (\ref{Ext1}) shows that up to conjugacy the $(p+1)$-tuple remains 
block-diagonal in first approximation w.r.t. $\varepsilon$. Hence, one 
can conjugate it by a matrix analytic in $\varepsilon$ to make the 
non-diagonal blocks zero in first approximation w.r.t. $\varepsilon$. 
In the same way 
one shows that the $(p+1)$-tuple is block-diagonal up to conjugacy of any 
order w.r.t. $\varepsilon$. The deformation being analytic, the 
$(p+1)$-tuple is block-diagonal up to conjugacy for $\varepsilon$ small 
enough and non-zero -- a contradiction.

$7^0$. If there exists $N\in {\cal Z}$ like in $1^0$, then one can conjugate 
the matrix $N$ and the matrices $M_j$ to the form
$N=\left( \begin{array}{ccc}0&0&I\\0&0&0\\0&0&0\end{array}\right)$, 
$M_j=\left( 
\begin{array}{ccc}P_j&U_j&V_j\\0&Q_j&W_j\\0&0&P_j\end{array}\right)$ 
where the middle row and column of blocks might be absent. If they are 
absent, then the monodromy 
group is a direct sum. Indeed, for the conjugacy classes $C_j'$ of the 
matrices $P_j$ one has $d(C_j')=d(C_j)/4$, hence, 
$d(C_1')+\ldots +d(C_{p+1}')=n^2/2$, see $5^0$. One has (\ref{Ext1}) with 
${\cal M}'={\cal M}''$ being the representation defined by the matrices $P_j$. 
Hence, the monodromy group is indeed a direct sum.

$8^0$. Suppose that the middle row and column of blocks are present. 
Lemma~\ref{Z} and 
conditions D) and E) imply that the multiplicities of the eigenvalues of the 
matrices $M_j$ 
for the diagonal blocks $P_j$ and $Q_j$ are proportional. Indeed, 
the blocks $S_j=P_j$ (the upper $P_j$) and 
$T_j=\left( \begin{array}{cc}P_j&U_j\\0&Q_j\end{array}\right)$ 
define invariant subspaces ${\cal X}_1$ and ${\cal X}_2$ of the 
monodromy group. Like in the case when $D\in {\cal Z}$, see $1^0$, and 
using the same notation one 
shows that equal eigenvalues of the matrices $S_j$ and $T_j$ are of 
proportional multiplicities. This implies that there holds 
(\ref{Ext1}), hence, the monodromy group is a direct sum of the groups defined 
by the blocks $S_j$ and $T_j$, i.e. after a simultaneous conjugation of the 
matrices $M_j$ one has $V_j=W_j=0$. Hence, there exists 
$D\in {\cal Z}$ like in $1^0$ which possibility 
is already rejected.

\subsection{Proof of Lemma \protect\ref{theMGisreducible}}

$1^0$. The monodromy group can 
be conjugated to a block upper-triangular form. The diagonal blocks define 
either irreducible or one-dimensional representations. The eigenvalues of 
each diagonal block $1\times 1$ 
satisfy the non-genericity relation $(\gamma _0)$ from Lemma~\ref{fivecond}. 

$2^0$. {\bf The lowest diagonal block is of size 1.} 

Indeed, set $M_j=\left( \begin{array}{cc}Q_j&\ast \\0&L_j\end{array}\right)$ 
where $L_j$ is the restriction of $M_j$ to the lowest diagonal block (say, of 
size $h$). 
Denote by $\Xi$, $\Theta$ the sets 
of eigenvalues $\sigma _{k,j}$, $j\leq p+1$ participating respectively 
in $(\gamma _B)$, $(\gamma _0)$ and by $\Xi '$, $\Theta '$ the sums of their 
respective eigenvalues $\lambda _{k,j}$. Hence, the set of eigenvalues 
of the blocks 
$L_1$, $\ldots$, $L_{p+1}$ is of the form $a\Xi +b\Theta$, $a\in {\bf N}$, 
$b\in {\bf Z}$. 

If $a>0$, $b\geq 0$, then condition $(\beta _h)$ is not fulfilled by the 
blocks 
$L_j$ (this condition is necessary because these blocks define an irreducible 
monodromy group of $h\times h$-matrices). Indeed, for $b=0$ it is not 
fulfilled because it is not fulfilled by the matrices $M_j$ and the 
multiplicities of equal eigenvalues of $M_j$ and $L_j$ are proportional.  
When increasing $h$, i.e. when increasing $b\in {\bf Z}$ while keeping 
$a$ fixed it is only the biggest 
multiplicity that increases and it is of an 
eigenvalue equal to 1. Hence, the sum of the quantities $r_j$ computed for 
the matrices $L_j$ remains the same while their size $h$ increases.  

On the other hand, one cannot have $b<0$ because in this case the sum 
of the eigenvalues $\lambda _{k,j}$ corresponding to the invariant solution 
subspace on which the monodromy group acts with the blocks $Q_j$ would be 
positive which contradicts Lemma~\ref{Z}. Indeed, the sum of these eigenvalues 
equals $\phi :=c\Xi '-b\Theta '+\Delta$ where $\Delta$ is the sum of some 
$n-h$ eigenvalues of the matrix $A_{p+2}$. We prove that $\phi >0$ like 
we prove that $\phi _1>0$ in $4^0$ of the proof of 
Lemma~\ref{fivecond}.

$3^0$. Denote by $\Pi$ the left upper $(n-1)\times (n-1)$-block. Conjugate 
it to make all non-zero rows of the restriction of 
the $(p+1)$-tuple $\tilde{M}$ of matrices $M_j-I$ to $\Pi$ linearly 
independent. After the conjugation some of the rows of the restriction of 
$\tilde{M}$ to $\Pi$ might be 0. In this case conjugate the matrices $M_j$ 
by one and the same permutation matrix which places the zero rows of 
$M_j-I$ in the last (say, $m$) positions (recall that the last row of 
$M_j-I$ is 0, see $2^0$, so $m\geq 1$). Notice that if the restriction 
to $\Pi$ of a row of  
$M_j-I$ is zero, then its last (i.e. $n$-th) position is 0 as well, 
otherwise $M_j$ is not diagonalizable. 
 
$4^0$. There remains to show that $m\geq n-n_1$. 
One has $M_j=\left( \begin{array}{cc}G_j&R_j\\0&I\end{array}\right)$, 
$I\in GL(m,{\bf C})$. 
Denote by $\tilde{G}$ the representation defined by the matrices $G_j$.
We regard the columns of 
the $(p+1)$-tuple of matrices $R_j$ as elements of the space 
${\cal F}(\tilde{G})$ (or just ${\cal F}$ for short) defined as follows.
Set $U^*=(U_1,\ldots ,U_{p+1})$. 
Set ${\cal D}=\{ U^*|U_j=(G_j-I)V_j,
V_j\in {\bf C}^m, \sum _{j=1}^{p+1}G_1\ldots G_{j-1}U_j=0\}$, 
${\cal E}=\{ U^*|U_j=(G_j-I)V,V\in {\bf C}^m\}$,  
${\cal F}={\cal D}/{\cal E}$.   

\begin{rem}
{\em If $R_j=(G_j-I)V$ with $V\in {\bf C}^m$ or with 
$V\in M_{m,n-m}$, then there holds} 
\begin{equation}\label{GR}
\sum _{j=1}^{p+1}G_1\ldots G_{j-1}R_j=0
\end{equation} 
{\em One has ${\cal E}\subset {\cal D}$. 
Equality (\ref{GR}) with $V\in M_{m,n-m}$ is condition 
(\ref{M_j}) restricted to the block $R$.}
\end{rem}

$5^0$. Each column of the $(p+1)$-tuple of matrices $R_j$ 
belongs to the linear space ${\cal D}$. 

{\bf The latter is of dimension 
$\theta =r_1+\ldots +r_{p+1}-(n-m)$.}

Indeed, the image of the linear operator 
$\tau _j:(.)\mapsto (G_j-I)(.)$ acting on ${\bf C}^{n-m}$ is of 
dimension $r_j$ (every column of $R_j$ 
belongs to the image of this operator, otherwise $M_j$ will not be 
diagonalizable). 
The $n-m$ 
linear equations resulting from (\ref{GR}) with 
$R_j=U_j=(G_j-I)V$, $V\in {\bf C}^m$ are linearly independent. 

Indeed, if they are not, then the images of all linear operators 
$\tau _j$ must be contained in a 
proper subspace of ${\bf C}^{n-m}$ (say, the one defined by the first 
$n-m-1$ vectors of its canonical basis). This means that all entries of the 
last rows of the matrices $G_j-I$ are 0. The matrices $M_j$ being 
diagonalizable, this implies that the entire $(n-m)$-th rows of $M_j-I$ are 
0. This contradicts the condition the first $n-m$ rows of the 
restriction to $\Pi$ of the $(p+1)$-tuple of matrices $M_j-I$ to be 
linearly independent, see $3^0$.

$6^0$. {\bf The space ${\cal F}$ is of codimension 
$n-m$ in ${\cal D}$, i.e. of dimension $\theta -2(n-m)$.}

Indeed, each vector-column $V$ 
belongs to ${\bf C}^{n-m}$ and the intersection ${\cal I}$ of the kernels of 
the operators $\tau _j$ is 
$\{ 0\}$, otherwise the matrices $M_j$ 
would have a non-trivial common centralizer -- if 
${\cal I}\neq \{ 0\}$, then after a change of the 
basis of ${\bf C}^{n-m}$ one can assume that a non-zero vector from 
${\cal I}$ equals 
$^t(1,0,\ldots ,0)$. Hence, the matrices $G_j$ are of the form 
$\left( \begin{array}{cc}1&\ast \\0&G_j^*\end{array}\right)$, 
$G_j^*\in GL(n-m-1,{\bf C})$, and one checks 
directly that $[M_j,E_{1,n}]=0$ for $E_{1,n}=\{ \delta _{i-1,n-j}\}$. 

$7^0$. The columns of the $(p+1)$-tuple of matrices $R_j$ (regarded 
as elements of ${\cal F}$) must be linearly independent, 
otherwise the monodromy group can be conjugated by a matrix 
$\left( \begin{array}{cc}I&\ast \\0&P\end{array}\right)$, 
$P\in GL(m,{\bf C})$, to a block-diagonal 
form in which the right lower blocks of $M_j$ are equal to $1$,  
the monodromy group is a direct sum and, hence, its centralizer is 
non-trivial -- a contradiction. 
This means that 
$\dim {\cal F}=\theta -2(n-m)=r_1+\ldots +r_{p+1}-2(n-m)\geq m$ 
which is equivalent to $m\geq n-n_1$; recall that 
$n_1=r_1+\ldots +r_{p+1}-n$. In the case of equality (and only in it) the 
columns of the $(p+1)$-tuple of matrices $R_j$ are a basis of the space 
${\cal F}$.~~~~$\Box$

\section{Case A)\protect\label{seccaseA}}

In this section we prove 

\begin{tm}\label{caseA}
The DSP is not solvable (hence, not weakly solvable, see Lemma~\ref{TC=Irr}) 
for quadruples of diagonalizable 
matrices $M_j$ each with MV equal to $(n/2,n/2)$ where $n\geq 4$ is even, the 
eigenvalues are relatively generic and $\xi$ is 
a non-primitive root of unity of order $n/2$. 
\end{tm}

\begin{rem}\label{caseArem}
In case A) for relatively generic eigenvalues 
there exist only block-diagonal quadruples of matrices $M_j$ with 
diagonal blocks $(n/l)\times (n/l)$. Their existence follows 
from \cite{Ko5}, Theorem 3. The non-existence of others follows from 
Theorem~\ref{caseA}.  
\end{rem}

The proof of the theorem consists of three steps. We assume 
that irreducible quadruples as 
described in the theorem exist. The first step is a preliminary deformation 
and conjugation of the quadruple which brings in some technical 
simplifications, the quadruple remaining irreducible and satisfying the 
conditions of the theorem, see the next subsection. At the second step we 
discuss the possible eigenvalues of the matrix $M_1M_2$ after the first step, 
see Subsection~\ref{evS}. At the third step we prove that 
the new quadruple must be reducible, see Subsection~\ref{endpr}.

\subsection{Preliminary conjugation and deformation}

Set 
$S=M_1M_2=(M_4)^{-1}(M_3)^{-1}$. Denote by $g_j,h_j$ the eigenvalues of $M_j$. 

\begin{lm}\label{but}
The triple $M_1,M_2,S^{-1}$ admits a  
conjugation to a block upper-triangular form with diagonal blocks of sizes 
only 1 or 2. The restriction of the triple to each 
diagonal block of size 2 is irreducible.
\end{lm} 

Indeed, suppose that the triple is in block upper-triangular form, its 
restrictions to each diagonal block being irreducible (in particular, the 
triple can be irreducible, i.e. with a single diagonal block). The restriction 
of $M_j$ to each diagonal block   
(say, of size $k$) is diagonalizable and has eigenvalues 
$g_j$ and $h_j$, of multiplicities $l^0$ and $k-l^0$. Hence, the 
conjugacy class of the restriction of $M_j$ to the block is of dimension 
$2l^0(k-l^0)\leq k^2/2$.

An irreducible 
triple with such blocks of $M_1$ and $M_2$ of size $k>1$ can exist only for  
$k=2$, in all other cases condition $(\alpha _k)$ does not hold. Indeed, 
the conjugacy class of the restriction of $S$ to the diagonal block is 
of dimension $\leq k^2-k$. Hence, the sum of the three 
dimensions is $\leq k^2/2+k^2/2+k^2-k=2k^2-k$ 
which is $<2k^2-2$ if $k>2$.~~$\Box$ \\ 

Give a more detailed description of the diagonal blocks of the triple 
$M_1,M_2,S^{-1}$ after the conjugation (in the form of 
lemmas; Lemmas~\ref{rank1}, \ref{semidirect} and \ref{Ext} are to be 
checked directly).  

\begin{lm}\label{rank1}
1) There are four possible representations defined by diagonal blocks 
of size 1 of the triple; we list them by indicating the couples of 
diagonal entries respectively of $M_1$ and $M_2$: 

\[ P~g_1,g_2~~~;~~~Q~h_1,h_2~~~;~~~R~g_1,h_2~~~;~~~U~h_1,g_2~~.\] 

2) Denote by $V$ and $W$ any two of these couples. For a given $V$ there 
exists 
a unique $W$ (denoted by $-V$) such that the corresponding diagonal entries 
of both $M_1$ and $M_2$ 
are different. One has $P=-Q$ and 
$R=-U$. 

3)One has dim~Ext$^1(V,W)=1$ if and only if $V=-W$. In the other cases one has 
dim~Ext$^1(V,W)=0$. 
\end{lm}

\begin{lm}\label{P=Q}
There are equally many diagonal blocks of type $V$ as there 
are of type $-V$. 
\end{lm}

Indeed, consider first the case when there are no blocks of size 2. 
Denote by $p'$, $q'$, $r'$ and $u'$ the number of blocks 
$P$, $Q$, $R$ and $U$. The multiplicities of the eigenvalues imply that 
$p'+r'=p'+u'=q'+u'=q'+r'=n/2$. Hence, $r'=u'$ and $p'=q'$.

If there are blocks of size 2, then each of them contains 
once each of the eigenvalues $g_1$, $g_2$, $h_1$, $h_2$ and the proof is 
finished in the same way as in the particular case considered above.~~$\Box$
 
\begin{lm}\label{rank2}
In an irreducible representation defined by a $2\times 2$-block the 
eigenvalues of $S$ can equal any 
couple $(\lambda ,\mu )$ (with $\lambda \mu =g_1h_1g_2h_2$) which is 
different from ($g_1g_2,h_1h_2$) and  
($g_1h_2,g_2h_1$).  
\end{lm}

Indeed, one can show (the easy computation is 
omitted) that if the eigenvalues of $S$ equal $g_1g_2$, $h_1h_2$ or 
$g_1h_2$, $g_2h_1$, then the triple is triangular up to conjugacy. 
On the other hand, if one fixes $M_1=$diag$(g_1,h_1)$ and varies $M_2$ 
within its conjugacy class, one can obtain any trace of the product $M_1M_2$. 
The determinant of the product being fixed, this means that $M_1M_2$ can 
belong to any non-scalar conjugacy class the product of whose eigenvalues 
equals $g_1h_1g_2h_2$. (The choice of the eigenvalues excludes the possibility 
$S$ to be scalar.)~~$\Box$

\begin{lm}\label{semidirect} 
The semi-direct sums defined by two diagonal blocks of size 1 are up to 
conjugacy of one of the types:  
$(M_1,M_2)=\left( \left( \begin{array}{cc}g_1&r\\s&h_1\end{array}\right) ,
\left( \begin{array}{cc}g_2&r'\\s'&h_2\end{array}\right) \right)~{\rm or}~
\left( \left( \begin{array}{cc}g_1&u\\m&h_1\end{array}\right) , 
\left( \begin{array}{cc}h_2&u'\\m'&g_2\end{array}\right) \right)$  
with either $r=r'=0$ or $s=s'=0$ but not both (resp. with either 
$u=u'=0$ or $m=m'=0$ but not both). Such semi-direct sums exist 
only for couples $(V,-V)$, see 2) and 3) from Lemma~\ref{rank1}. The 
centralizers of these semi-direct sums are trivial. 
\end{lm}

Denote by $\Phi$, $\Psi$ respectively an irreducible representation 
of rank 2 defined by a diagonal 
block of the triple $M_1,M_2,S^{-1}$ and a representation which is either 
irreducible and non-equivalent to $\Phi$ or one-dimensional (i.e. of type 
$P$, $Q$, $R$ or $U$, see 
Lemma~\ref{rank1}) or a semi-direct sum of two one-dimensional ones $(V,-V)$, 
see Lemmas~\ref{rank1} and \ref{semidirect}.

\begin{lm}\label{Ext}
One has dim~Ext$^1(\Phi ,\Psi )=$dim~Ext$^1(\Psi ,\Phi )=0$.
\end{lm}

\begin{defi} 
{\em We say that the triple $M_1,M_2,S^{-1}$ or $M_3,M_4,S$ is in a 
{\em special form} if it is block-diagonal,  
each diagonal block $B_{\mu}$ being itself block upper-triangular, its 
diagonal blocks being of equal size which is either 1 or 2. In the case of 
size 2 all diagonal blocks of each block $B_{\mu}$ define equivalent 
representations. 
In the case of size 1 the block $B_{\mu}$ is of size 2 and defines a 
semi-direct sum, see Lemma~\ref{semidirect}. Thus a triple in special form 
is block upper-triangular with diagonal blocks of size 2 defining either 
irreducible representations or semi-direct sums like in 
Lemma~\ref{semidirect}.}
\end{defi} 

\begin{lm}\label{specialform}
One can deform the matrices $M_j$ within their conjugacy classes 
(without changing the matrix $S$) 
so that after the deformation each  
of the triples $M_1,M_2,S^{-1}$ and $M_3,M_4,S$ after a suitable conjugation 
is in special form. The two conjugations are, in general, different.
\end{lm} 

The lemma is proved in Subsection \ref{proofspecialform}.

\subsection{The possible eigenvalues of the matrix $S$\protect\label{evS}} 

The eigenvalues of the matrix $S$ (even when they are distinct) 
must satisfy certain equalities -- for every diagonal 
block of size 2 (irreducible or not) 
of the triple $M_1,M_2,S^{-1}$ (resp. $M_3,M_4,S$) the eigenvalues 
$\lambda ,\mu$ of $S$  must satisfy the 
condition $g_1h_1g_2h_2\lambda ^{-1}\mu ^{-1}=1$ 
(resp. $g_3h_3g_4h_4\lambda \mu =1$). 

In what follows we denote the 
eigenvalues of $S$ by $s_i$.  
Let the triple $M_1,M_2,S^{-1}$ (resp. $M_3,M_4,S$) be in special form. 
For each eigenvalue $s_i$ denote by 
$t(s_i)$ (resp. by $u(s_i)$) the eigenvalue of $S$ in the same diagonal 
$2\times 2$-block of the 
triple with 
$s_i$. Note  
that $t(t(s_i))=s_i=u(u(s_i))$. One has $t(s_i)=u(s_i)$ if and only if 
$\xi =1$ (and this holds for all $i=1,\ldots ,n$).\\

Set $i_1=1$. For the eigenvalue $s_1=s_{i_1}$ find 
$s_{i_2}\stackrel{{\rm def}}{=}t(s_{i_1})$, then find 
$s_{i_3}\stackrel{{\rm def}}{=}u(s_{i_2})$, then 
$s_{i_4}\stackrel{{\rm def}}{=}t(s_{i_3})$, then again 
$s_{i_5}\stackrel{{\rm def}}{=}u(s_{i_4})$ etc. Thus one has 
$s_{i_{\nu +1}}=t(s_{i_{\nu}})$ for $\nu$ odd (hence, 
$t(s_{i_{\nu +1}})=t(t(s_{i_{\nu}}))=s_{i_{\nu}}$) and 
$s_{i_{\nu +1}}=u(s_{i_{\nu}})$ for $\nu$ even (hence, 
$u(s_{i_{\nu +1}})=u(u(s_{i_{\nu}}))=s_{i_{\nu}}$. 

Denote by $m$ the least value 
of $\alpha$ for which one has $i_{\alpha}=1$. It is clear that $m-1$ is even.

\begin{lm}\label{xixi-1}
For $\nu$ odd one has $s_{i_{\nu +1}}=\xi s_{i_{\nu -1}}$, for 
$\nu$ even one has $s_{i_{\nu +1}}=\xi ^{-1}s_{i_{\nu -1}}$.
\end{lm}

Indeed, there holds $g_1h_1g_2h_2s_{i_{\nu}}^{-1}(t(s_{i_{\nu }}))^{-1}=
g_3h_3g_4h_4s_{i_{\nu}}u(s_{i_{\nu }})=1~{\rm and~}
\prod _{j=1}^4g_jh_j=\xi$.  
Hence, $\xi ^{-1}t(s_{i_{\nu}})=u(s_{i_{\nu}})$. 
For $\nu$ odd this yields 
$\xi ^{-1}s_{i_{\nu +1}}=\xi ^{-1}t(s_{i_{\nu}})=u(s_{i_{\nu}})=
u(u(s_{i_{\nu -1}}))=s_{i_{\nu -1}}$, for $\nu$ even in the same way it gives 
$\xi ^{-1}s_{i_{\nu -1}}=s_{i_{\nu +1}}$.~~$\Box$ 

\begin{lm}\label{mn}
One has $m-1<n/2$ and $m-1$ divides $n/2$.
\end{lm}

{\em Proof:} Recall that $\xi =\exp (2k\pi i/(n/2))=\exp (4k\pi i/n)$ (see 
Subsection~\ref{GEPsi}). 
If $k=0$, i.e. $\xi =1$, 
then $s_3=s_1$, i.e. $m-1=1$, and the statement holds.  

Let $k\neq 0$. Then $s_{1+(m-1)}=(\xi )^{-m+1}s_1=s_1$ (Lemma~\ref{xixi-1}).  
Hence, $(\xi )^{-m+1}=1$, i.e.   
$4k(m-1)=2nl$ ($l$ is defined in Subsection~\ref{GEPsi}), i.e. 
$k(m-1)=(n/2)l$. The minimality of 
$m$ (hence, of $m-1$ as 
well) implies that $m-1$ and $l$ are relatively prime, i.e. $m-1$ divides 
$n/2$. The non-primitivity 
of $\xi$ implies $k>1$. Hence, $m-1<n/2$.~~$\Box$ 

\begin{rem}\label{partition}
{\em Lemma \ref{mn} implies that the set of eigenvalues of $S$ can be 
partitioned into $n/(2m-2)$ sets ${\cal N}_1$, $\ldots$, 
${\cal N}_{n/(2m-2)}$  
each consisting of $(2m-2)$ eigenvalues (denoted again by $s_i$)  
with the properties $s_{2k+2}=\xi s_{2k}$, $s_{2k+1}=\xi ^{-1}s_{2k-1}$,   
$s_{2k-1}s_{2k}=g_1h_1g_2h_2$ and 
$s_{2k}^{-1}s_{2k+1}^{-1}=g_3h_3g_4h_4$. If some of the sets ${\cal N}_i$ 
are identical, then we define their {\em multiplicities} 
in a natural way. Two non-identical sets ${\cal N}_i$ have no eigenvalue in 
common. In what follows we change the indexation -- equal (different) 
indices indicate identical (different) sets ${\cal N}_i$.} 
\end{rem}

\subsection{End of the proof of Theorem \protect\ref{caseA}
\protect\label{endpr}} 

{\bf Case 1)} {\em The matrix $S$ has  
at least two different sets ${\cal N}_i$.}
 
Then the upper-triangular form of 
the triple $M_1,M_2,S$ is in addition block-diagonal, the restrictions 
of the matrix $S$ to two different diagonal blocks having no eigenvalue in 
common. Indeed, it suffices to rearrange the blocks $B_{\mu}$ from the 
special form putting first all the blocks  $B_{\mu}$ with eigenvalues of 
$S$ from 
${\cal N}_1$ (repeated with its multiplicity -- 
this defines the diagonal block $R_1$), then all blocks with 
eigenvalues of $S$ from ${\cal N}_2$ (this defines the diagonal block $R_2$) 
etc. The size of the block $R_i$ equals $l_i$ times the number of 
eigenvalues from ${\cal N}_i$, $l_i\in {\bf N}^*$.

The triple $M_3,M_4,S$ admits a conjugation to the same block-diagonal form. 
Hence, if the triple $M_1,M_2,S^{-1}$ is block-diagonal (with diagonal blocks 
$R_i$), to give the same form of the triple $M_3,M_4,S$ one has to use 
as conjugation matrix one commuting with $S$, hence, a 
block-diagonal one with diagonal blocks of the sizes of the blocks $R_i$. 
Hence, both triples are simultaneously block-diagonal, i.e. the 
quadruple 
$M_1$, $M_2$, $M_3$, $M_4$ is block -diagonal, i. e. reducible.\\ 

{\bf Case 2)} {\em There is a single set ${\cal N}_1$ 
repeated $n/(2m-2)$ times.}
In this case one can deform the matrices $M_j$, $j=1,2$, so that the matrix 
$S$ have at least two different sets ${\cal N}_i$ of eigenvalues. 

\begin{defi}
{\em We say that a matrix is in $s$-block-diagonal (resp. in $s$-block 
upper-triangular) form if it is block-diagonal (resp. block 
upper-triangular) with diagonal blocks all of size $s$.}
\end{defi}

Set $\mu =n/(2m-2)$. Conjugate the triple $M_1,M_2,S$ to 
a $(2m-2)$-block upper-triangular form where the diagonal blocks of the 
matrix $S$ are with eigenvalues from ${\cal N}_1$: 

\[ M_j=\left( \begin{array}{cccc}M_j'&H_{j;1,2}&\ldots &H_{j;1,\mu}\\
0&M_j'&\ldots &H_{j;2,\mu}\\ \vdots &\vdots &\ddots &\vdots \\0&0&\ldots &M_j'
\end{array}\right) ~~,~~j=1,2~~,~~
S=\left( \begin{array}{cccc}T&Q_{1,2}&\ldots &Q_{1,\mu}\\
0&T&\ldots &Q_{2,\mu}\\ \vdots &\vdots &\ddots &\vdots \\0&0&\ldots &T
\end{array}\right) ~~.\]
We assume that the blocks $M_j'$ and $T$ are 2-block-diagonal. 

Deform analytically the left upper blocks of size $2m-2$ of the matrices 
$M_1$, $M_2$ and $S$ so that they remain 2-block-diagonal and the eigenvalues 
of $S$ change to new ones, forming again a set of $2m-2$ eigenvalues like 
in Remark~\ref{partition} but different from ${\cal N}_1$. To this end one 
can keep the matrix $M_1$ the same and vary the left upper block of the 
matrix $M_2$; see Lemma~\ref{rank2}. 
This block will become $M_2'+\varepsilon U$, 
$\varepsilon \in ({\bf C},0)$, $U\in gl(2m-2,{\bf C})$, and the one of 
$S$ will equal 
$M_1(M_2'+\varepsilon U)$. The other blocks of $M_1$, $M_2$ and $S$ do not 
change. 

One can deform in a similar way the triple of matrices $M_3^{-1}$, $M_4^{-1}$ 
$S$ (requiring the deformation of $S$ to be the same in both triples). For 
$\varepsilon \neq 0$ small enough the quadruple of matrices remains 
irreducible. However, there are already two different sets ${\cal N}_i$ of 
eigenvalues of $S$, so 
we are in Case 1) and the quadruple is block-diagonal. Hence, the initial 
quadruple is also reducible.

\subsection{Proof of Lemma 
\protect\ref{specialform}\protect\label{proofspecialform}}

\begin{nota}
{\em Assume that the triple $M_1,M_2,S$ satisfies the conclusion of 
Lemma~\ref{but}. Block-decompose each matrix 
from $gl(n,{\bf C})$ the sizes of the diagonal blocks being the same as the 
ones of the triple $M_1,M_2,S$. Denote the block of this decomposition in 
the $i$-th row and 
$k$-th column of blocks by $([i,k])$. By $(i,k)$ we denote the matrix entry in 
the $i$-th row and $k$-th column.}
\end{nota}

$1^0$. {\em Up to conjugacy the triple $M_1,M_2,S$ is block-diagonal, 
with two diagonal blocks ($T$ and $Y$) which are block upper-triangular, 
their diagonal 
blocks being respectively of size 1 and 2, the latter defining irreducible 
representations.}

Indeed, whenever a block $([i,i+1])$ of the triple $M_1,M_2,S$ 
is of size $1\times 2$ or $2\times 1$, it can be made equal to 0 by a 
simultaneous conjugation of the triple with a matrix of the form $I+R$ 
where only the block $([i,i+1])$ of $R$ is non-zero. This follows from 
Lemma~\ref{Ext}. After this in the same way one annihilates all blocks 
$([i,i+2])$ of size $1\times 2$ or $2\times 1$, then all blocks $([i,i+3])$ 
of these sizes etc. Then one rearranges the diagonal blocks putting the ones 
of size 1 first and the ones of size 2 next. This gives the claimed form.\\

$2^0$. {\em The block $Y$ after conjugation becomes block-diagonal, its 
diagonal blocks $B_{\mu}$ being block upper-triangular, their diagonal blocks 
being of size 2. The diagonal blocks of one and the same (resp. of different) 
diagonal blocks $B_{\mu}$ define equivalent (resp. non-equivalent) 
representations.}

This is proved by analogy with $1^0$, making use of Lemma \ref{Ext}.\\  

$3^0$. Denote by $V_1$, $\ldots$, $V_n$ the diagonal blocks of $T$. 

{\em One can conjugate the 
triple $M_1,M_2,S$ by an upper-triangular matrix so that after the 
conjugation only these blocks $([i,j])$, $i<j$, remain possibly non-zero 
for which $V_i=-V_j$.}

This is proved like $1^0$ and $2^0$, making use of 2) and 3) 
of Lemma~\ref{rank1}.\\   

$4^0$. {\em After a conjugation and deformation the block $T$ of the triple 
$M_1,M_2,S$ becomes block-diagonal, with upper-triangular diagonal blocks of 
size 2 defining semi-direct sums, see Lemma~\ref{semidirect}.} 

The proof of this statement occupies $4^0$ -- $5^0$. It completes the proof 
of the lemma.

A conjugation of the triple $M_1,M_2,S$ with a permutation 
matrix 
places the set of blocks $P$ and $Q$ first and the set of blocks $R$ and $U$ 
last on the diagonal; the triple remains block upper-triangular, in addition 
it is block-diagonal, the sizes of the diagonal blocks equal respectively 
$\sharp P+\sharp Q$ and $\sharp R+\sharp U$ (one of these sizes can be 0).  

It suffices to consider the case when only, say, blocks $P$ and $Q$ 
are present, in the general case the reasoning is the same. Observe first 
that the blocks $P$ and $Q$ can be situated on the diagonal in any possible 
order. 

The eigenvalues of the restrictions of $S$ to the blocks $P$ and $Q$ being 
different, one can conjugate the triple with an upper-triangular matrix 
to make $S$ diagonal. Moreover, all blocks $([i,j])$, $i<j$, with $V_i=V_j$ 
are 0, otherwise 
at least one of the matrices $M_1$, $M_2$ will not be diagonal.\\  

$5^0$. Consider first the case when the triple after this conjugation becomes 
diagonal. Rearrange the blocks in alternating order -- $P$, $Q$, $P$, $Q$, 
$\ldots$. Make non-zero the entries $(1,2)$, $(3,4)$, $(5,6)$ etc. of the 
matrices $M_j$ without changing the matrix $S$. With the notation from 
Lemma~\ref{semidirect} this amounts to choosing $s=s'=0$, $r\neq 0$, 
$r'=-rh_2g_1^{-1}$ (look at the first couple $(M_1,M_2)$ from the lemma). 
This gives the necessary block-diagonal form of the block $T$. The 
representations $P$ and $Q$ being non-equivalent, the centralizers of the 
diagonal blocks are trivial.

Suppose now that the triple is not diagonalizable and that $V_1=P$ 
(the case $V_1=Q$ is considered by analogy).
Denote by $i_1<\ldots <i_h$ the indices $i$ for which $V_i=Q$. Denote by 
$m$ the {\em smallest} $i_{\nu }$ for which at least one of the entries 
$(k,m)$ of $M_1$ and $M_2$ is non-zero, $k<m$; by $3^0$, $k$ is not among the 
indices $i_{\nu }$. Denote the {\em greatest} such value of $k$ by $k_0$. 
Hence, all entries $(i,k_0)$ ($i<k_0$) and $(k_0,\mu )$ ($\mu <m$) 
of $M_1$ and $M_2$ are 0, otherwise these 
matrices will not be diagonalizable. 

One can annihilate all entries $(k',m)$ of $M_j$ where  
$k'<k_0$ by consecutively conjugating the triple $M_1,M_2,S$ by matrices of 
the form $I+gE_{k',k_0}$. Note that the values of $k'$ are not among the 
indices $i_{\nu }$. In a similar way one annihilates all entries $(k_0,k'')$ 
of $M_j$ with $k''>m$ by consecutive conjugations with matrices of the form 
$I+gE_{m,k''}$. 

Hence, it is possible to conjugate the triple by a permutation matrix 
putting the $k_0$-th and $m$-th rows and columns first and preserving its 
upper-triangular form; in addition, the triple will be block-diagonal with 
first diagonal block of size 2 (which is upper-triangular  
non-diagonal and with trivial centralizer). After this one continues in the 
same way with the lower block. 
In the end the block $T$ will become upper-triangular and block-diagonal, 
with diagonal blocks of size 2 each of which is triangular non-diagonal with 
trivial centralizer. 

\section{Case C)\protect\label{seccaseC}}

\begin{lm}\label{helpful}
If $\kappa =0$ and if the DSP is solvable for a $(p+1)$-tuple of conjugacy 
classes $C_j$ with relatively generic eigenvalues 
defining the diagonal JNFs $J_j^n$, then the DSP is solvable for any 
$(p+1)$-tuple of JNFs ${J_j'}^n$ and for any relatively generic eigenvalues 
with the same value of $\xi$ 
where for each $j$ the JNFs $J_j^n$ and ${J_j'}^n$ correspond to one another 
or are the same. 
\end{lm}

The lemma is proved at the end of the subsection.

Assume that there exist irreducible triples of diagonalizable matrices $M_j$ 
such that $M_1M_2M_3=I$, the PMV of the eigenvalues of the matrices being 
equal to $(d,d,d,d)$, $(d,d,d,d)$, $(2d,2d)$. Denote by $\sigma _{k,j}$ 
the eigenvalues of $M_j$ where $k=1,2,3,4$ if $j=1$ or 2 and $k=1,2$ if 
$j=3$. 

One can choose the eigenvalues 
of $M_1$ and $M_2$ such that $\sigma _{1,j}=-\sigma _{2,j}$ and 
$\sigma _{3,j}=-\sigma _{4,j}$, $j=1,2$, see Lemma~\ref{helpful}. 
Hence, the MVs of the eigenvalues 
of the matrices $(M_1)^2$ and $(M_2)^2$ are of the form $(2d,2d)$. 
Set $A=M_1M_2=(M_3)^{-1}$, $B=M_2M_1$. The matrix $B$ 
is conjugate to $(M_3)^{-1}$ (because $B=M_2(M_3)^{-1}(M_2)^{-1}$).  
One has $AB=M_1(M_2)^2M_1$, 
hence, $AB=(M_1)^2(M_1)^{-1}(M_2)^2M_1$. Set 

\[ L_1=A=M_1M_2~~,~~L_2=B=M_2M_1~~,~~L_3=(M_1)^{-1}(M_2)^{-2}M_1~~,~~
L_4=(M_1)^{-2}~.\]
One has $L_1L_2L_3L_4=I$. The matrices $L_j$ are diagonalizable, their 
MVs equal $(2d,2d)$ and by Case A) they 
define a block-diagonal algebra ${\cal C}$ with $2k$ blocks $2s\times 2s$. 
Hence, dim${\cal C}\leq 8ks^2$. 

The algebra ${\cal C}$ contains also the matrices $(L_j)^{-1}$. 
Hence, it contains the matrices $(M_1)^2=(L_4)^{-1}$, $M_1M_2=L_1$,  
$M_2M_1=L_2$ and $(M_2)^2=(M_2M_1)(L_3)^{-1}(M_2M_1)^{-1}$. 

Every matrix from the algebra ${\cal D}$ 
generated by $M_1$ and $M_2$ is 
of the form $K+M_1L+M_2S$ with $K,L,S\in {\cal C}$. Hence, 
dim${\cal D}\leq 3$dim${\cal C}<n^2=$dim$gl(n,{\bf C})$. By the 
Burnside theorem, the matrix algebra ${\cal D}$ 
is reducible. 

{\bf Proof of Lemma~\ref{helpful}:}
$1^0$. Suppose that the DSP is not solvable for the JNFs ${J_j'}^n$ and for 
some relatively generic but not generic eigenvalues. Prove that 
then it is not solvable for 
the JNFs $J_j^n$ and for any such eigenvalues. 
Note first that the JNFs $J_j^n$ 
and ${J_j'}^n$ satisfy the conditions of Theorem~\ref{generic}, see 
Corollary~\ref{corpsicorr}. 

$2^0$. An irreducible $(p+1)$-tuple ${\cal H}$ of matrices $M_j$ with JNFs 
$J_j^n$ can be realized by a Fuchsian system with diagonalizable 
matrices-residua $A_j$ such that $J(A_j)=J(M_j)$ for $j\leq p+1$ and with 
an additional apparent singularity, 
see Subsection~\ref{pritself} with the definition of the sets ${\cal G}_i$, 
the maps $\chi _i$ and ${\cal M}$. One can choose $i$ such that 
$\chi _i({\cal G}_i)$ is dense in ${\cal M}$.

$3^0$. Vary the eigenvalues of the matrices $A_j$ within the set ${\cal G}_i$ 
without changing their 
JNFs. For suitable eigenvalues (in general, with integer differences 
between some of them; see $5^0$) 
one obtains as monodromy group ${\cal H}'$ 
of the Fuchsian system one in which either 
$J(M_j)={J_j'}^n$ or $J(M_j)$ is {\em subordinate} to ${J_j'}^n$, i.e. 
the multiplicities of the eigenvalues are the same and 
for each eigenvalue $\lambda$ and for each 
$s\in {\bf N}$ rk$(M_j-\lambda I)^s$ is the same or smaller than 
should be, see the 
details in \cite{Ko4}. One can assume that 
the eigenvalues of the matrices $M_j$ are relatively generic. 
Such a monodromy 
group cannot be irreducible (otherwise one could deform it using the basic 
technical tool into a nearby one with the same eigenvalues 
and with $J(M_j)={J_j'}^n$ for all $j$; such irreducible monodromy groups do 
not exist by assumption). 

$4^0$. The monodromy group ${\cal H}'$ can be analytically deformed into the 
monodromy group ${\cal H}$ because both are obtained from the Fuchsian 
system for different eigenvalues of the matrices-residua. However, 
${\cal H}'$ cannot be analytically deformed into a nearby irreducible 
monodromy group with JNFs as in ${\cal H}$. 

Indeed, if for all $j$ one has $J(M_j)={J_j'}^n$ in 
${\cal H}'$, then the monodromy group 
${\cal H}'$ must be block-diagonal with diagonal blocks of equal size 
and for the representations 
$\Phi _1$, $\Phi _2$ defined by two diagonal blocks one has 
Ext$^1(\Phi _1,\Phi _2)\leq 0$ with equality if and only if $\Phi _1$, 
$\Phi _2$ are not equivalent. The last inequality holds also if for 
some $j$ $J(M_j)$ 
is subordinate to ${J_j'}^n$. After this one applies the reasoning from 
$5^0$ -- $8^0$ of the proof of Lemma~\ref{fivecond}. 

$5^0$. It is explained in \cite{Ko4} how to choose the eigenvalues from $3^0$ 
to obtain the monodromy group ${\cal H}'$ with $J(M_j)$ equal or 
subordinate to ${J_j'}^n$. Their possible choice is not unique -- if one 
adds to equal eigenvalues of the matrices $A_j$ equal integers the sum of all 
added integers (taking into account the multiplicities) being 0, then one 
obtains a new possible such set of eigenvalues; different eigenvalues of a 
given matrix $A_j$ must 
remain such and if two eigenvalues of a given matrix $A_j$ differ by a 
non-zero integer, then the order of their real parts must be preserved. 

From all these a priori possible choices there is at least one which is really 
possible, i.e. for which there exists such a point from ${\cal G}_i$. Indeed, 
${\cal G}_i$ is constructible and its projection on the set of eigenvalues 
${\cal W}$ must be dense in ${\cal W}$, see 
Subsection~\ref{pritself}.~~~~$\Box$

\section{Case B)\protect\label{seccaseB}}
  
\begin{defi}
{\em A {\em special triple} is an irreducible triple of 
matrices $M_j$ such that 
$M_1M_2M_3=I$, $M_1-I$ and $M_2-I$ being conjugate to nilpotent Jordan 
matrices consisting each of $n/3$ Jordan blocks of size 3, $M_3$ being 
diagonalizable, with three eigenvalues each of multiplicity $n/3$. The 
eigenvalues are presumed to be relatively generic but not generic.}
\end{defi} 

In the present subsection we prove that special triples do not exist. By 
Lemma~\ref{helpful}, there exist no irreducible triples from Case B).

\begin{lm}\label{dirsum}
Suppose that there exist special triples. Then there exist special triples 
satisfying the conditions

{\em i)} Im$(M_j-I)\cap {\rm Ker}(M_{2-j}-I)=\{ 0\}$, 
$j=1,2$

{\em ii)} ${\bf C}^n={\rm Ker}(M_1-I)\oplus {\rm Ker}(M_2-I)\oplus 
({\rm Im}(M_1-I)\cap {\rm Im}(M_2-I))$.
\end{lm}

\begin{cor}\label{dirsumcor}
If there exist special triples, then there exist special triples in which the 
matrices $M_1-I$, $M_2-I$ are of the form 
$M_1-I=\left( \begin{array}{ccc}0&0&0\\P&0&0\\Q&R&0\end{array}\right)$, 
$M_2-I=\left( \begin{array}{ccc}0&I&V\\0&0&I\\0&0&0\end{array}\right)$ 
in which all blocks are $(n/3)\times (n/3)$, the matrices $P$ and $R$ being 
non-degenerate.
\end{cor}

The lemmas and the corollary from this section are proved at its end. 
Let the matrices $M_j$ be like in Corollary~\ref{dirsumcor}. Consider 
the matrices 

\[ N_1=\left( \begin{array}{ccc}I&0&V-I\\0&I&I\\0&0&I\end{array}\right) ,~
N_2=\left( \begin{array}{ccc}I&0&0\\P&I&0\\Q&0&I\end{array}\right) ,~ 
G=\left( \begin{array}{ccc}0&I&0\\0&0&0\\0&0&0\end{array}\right),~
H=\left( \begin{array}{ccc}0&0&0\\0&0&0\\0&R&0\end{array}\right) .\] 
Hence, each of the matrices $N_1-I$, $N_2-I$, $G+H$, $G$ and $H$ is 
nilpotent and conjugate to a Jordan matrix 
consisting of $n/3$ blocks of size 2 and of $n/3$ blocks of size 1. One has 
(to be checked directly) 

\begin{equation}\label{relations1a} 
N_2(I+H)=M_1~~,~~{\rm i.e.}~M_1-I=(N_2-I)(I+H)+(G+H)-G
\end{equation}
\begin{equation}\label{relations1b}
(I+G)N_1=M_2~~,~~{\rm i.e.}~M_2-I=(I+G)(N_1-I)+G
\end{equation}
\begin{equation}\label{relations2} 
GH=G^2=H^2=HG=0~~,~~
(N_1-I)G=0~~,~~(N_2-I)H=0
\end{equation}
Hence, $N_2(I+G+H)N_1=M_1M_2$. 
Denote by ${\cal A}$ the matrix algebra generated by the matrices $N_1-I$, 
$G+H$ and $N_2-I$, by ${\cal B}$ the one generated by $M_1$ and $M_2$. 

\begin{lm}\label{dimA}
The matrix algebra ${\cal A}$ is reducible and dim${\cal A}\leq n^2/2$.
\end{lm}

One has 
${\cal B}={\cal A}+G{\cal A}+{\cal A}G+G{\cal A}G+
G{\cal A}G{\cal A}+{\cal A}G{\cal A}G+\ldots ~(*)$.
Indeed, every product of the matrices $M_1-I$ and $M_2-I$ 
(in any order and  
quantity) is representable as a linear combination of such products of the 
matrices $N_1-I$, $N_2-I$, $G+H$ and $G$, 
see (\ref{relations1a}), (\ref{relations1b}) and (\ref{relations2}). 

On the other hand, one has ${\cal A}G\subset {\cal A}$. Indeed, denote by  
$Y$ a product of the matrices $N_1-I$, $N_2-I$, $G+H$ (in any quantity and 
order). If its right most factor is $N_1-I$ or $G+H$, then by 
(\ref{relations2}) one has $YG=0$. If it is $N_2-I$, then  
$YG=Y(G+H)-YH=Y(G+H)\in {\cal A}$. 

This together with $(*)$ implies that ${\cal B}={\cal A}+G{\cal A}~(**)$.  
Suppose that the couple of matrices 
$M_1$, $M_2$ is irreducible. Then by the Burnside theorem 
the algebra ${\cal B}$ equals $gl(n,{\bf C})$, i.e. dim${\cal B}=n^2$. The 
restriction of each matrix from ${\cal B}$ to the last $2n/3$ rows is the 
restriction to them of a matrix from ${\cal A}$, see $(**)$. This 
means that dim${\cal A}\geq 2n^2/3$ which contradicts Lemma~\ref{dimA}. Hence, 
special triples do not exist.

{\bf Proof of Lemma \ref{dirsum}:}
$1^0$. Recall that the three conjugacy classes $C_j$ of the matrices $M_j$ 
belong to $SL(n,{\bf C})$. 
Denote by ${\cal U}$ the variety of irreducible  
representations (i.e. triples $(M_1,M_2,M_3)$ defined up to conjugacy) 
where $M_j\in C_j\subset SL(n,{\bf C})$, 
$M_1M_2M_3=I$. 

Find dim${\cal U}$. One has to consider the cartesian product 
$C_1\times C_2\subset (SL(n,{\bf C})\times SL(n,{\bf C}))$. 
The algebraic variety ${\cal V}\subset (SL(n,{\bf C}))^2$ of irreducible 
couples of 
matrices $M_1$, $M_2$ such that $M_1\in C_1$, $M_2\in C_2$ and 
$(M_1M_2)^{-1}\in C_3$  
is the projection in $C_1\times C_2$ of the intersection of 
the two varieties in $C_1\times C_2\times SL(n,{\bf C})$: the cartesian 
product $C_1\times C_2\times C_3$ and the graph of the 
mapping 
$(C_1\times C_2)\ni (M_1, M_2)
\mapsto M_3=M_2^{-1}M_1^{-1}\in SL(n,{\bf C})$. This intersection is 
transversal which implies the smoothness of the 
variety ${\cal V}$ (this can be proved by analogy with 1) of Theorem 2.2 
from \cite{Ko2}). Thus 
${\rm dim}\, {\cal V}=(\sum _{j=1}^2{\rm dim}\, C_j)-[(n^2-1)-
{\rm dim}\, C_3]$
(here $(n^2-1)-{\rm dim}\, C_3=$codim$_{SL(n,{\bf C})}C_3$). 
Hence, dim${\cal V}=$dim$\, C_1+$dim$\, C_2+$dim$\, C_3-n^2+1$. 

$2^0$. In order to obtain dim${\cal U}$ from dim${\cal V}$ 
one has to factor out the possibility to 
conjugate the triple $(M_1,M_2,M_3)$ with matrices from 
$SL(n,{\bf C})$. No non-scalar such matrix commutes with all the 
matrices $(M_1,M_2,M_3)$ due to the irreducibility of the 
triple and to Schur's lemma. Thus 
${\rm dim}\, {\cal U}=
{\rm dim}\, {\cal V}-{\rm dim}\, SL(n,{\bf C})=
\sum _{j=1}^3{\rm dim}\, C_j-2n^2+2=2$.

$3^0$. The subvariety ${\cal U}'\subset {\cal U}$ on which one has 
dim~(Ker$(M_j-I)\cap$Im$(M_{2-j}-I))>0$ for $j=1$, 2 is of positive 
codimension in ${\cal U}$. Indeed, its dimension is computed like the one 
of ${\cal U}$, 
by replacing the cartesian product $C_1\times C_2$ by its subvariety on which 
one has dim~(Ker$(M_j-I)\cap$Im$(M_{2-j}-I))>0$ for $j=1$, 2. This  
subvariety is of positive codimension. Hence, the condition 
dim~(Ker$(M_j-I)\cap$Im$(M_{2-j}-I))>0$ for $j=1$, 2 cannot hold for all 
points from ${\cal U}$. 

Condition {\em ii)} follows from condition {\em i)}. ~~~~$\Box$  

{\bf Proof of Corollary \ref{dirsumcor}:}
$1^0$. One has dim~Ker$(M_1-I)=$dim~Ker$(M_2-I)=n/3$. 
Condition {\em ii)} of Lemma~\ref{dirsum} implies that 
dim~(Im$(M_1-I)\cap {\rm Im}(M_2-I))=n/3$; recall that 
Ker$(M_j-I)\subset$Im$(M_j-I)$, $j=1,2$. Choose a basis of ${\bf C}^n$ 
such that the first $n/3$ vectors are a basis of ${\rm Ker}(M_2-I)$, the next 
$n/3$ vectors are a basis of ${\rm Im}(M_1-I)\cap {\rm Im}(M_2-I)$ and the 
last $n/3$ vectors are a basis of ${\rm Ker}(M_1-I)$. Hence, in this basis the 
matrices of $M_1-I$, $M_2-I$ look like this: 
$M_1-I=\left( \begin{array}{ccc}0&0&0\\P'&T&0\\Q'&R'&0\end{array}\right)~~,~~
 M_2-I=\left( \begin{array}{ccc}0&W&V'\\0&U&Y\\0&0&0\end{array}\right)$ 
(all blocks are $(n/3)\times (n/3)$).

$2^0$. One has $(M_2-I)^3=\left( 
\begin{array}{ccc}0&WU^2&WUY\\0&U^3&U^2Y\\0&0&0\end{array}\right) =0$. The 
rank of the matrix $\left( \begin{array}{c}W\\U\end{array}\right)$ equals 
$n/3$ because rk$(M_2-I)=2n/3$. Therefore the equalities 
$\left( \begin{array}{c}WU^2\\U^3\end{array}\right)=
\left( \begin{array}{c}0\\0\end{array}\right)$ and   
$\left( \begin{array}{c}WUY\\U^2Y\end{array}\right)=
\left( \begin{array}{c}0\\0\end{array}\right)$ 
imply respectively $U^2=0$ and $UY=0$. It follows from  rk$(M_2-I)=2n/3$ 
that rk$(U~Y)=n/3$. Hence, the equality $(U^2~UY)=(0~0)$ implies $U=0$. 

$3^0$. In the same way one proves that $T=0$. A simultaneous 
conjugation of $M_1-I$ and $M_2-I$ with the matrix 
$\left( \begin{array}{ccc}WY&0&0\\0&Y&0\\0&0&I\end{array}\right)$ 
brings them to the desired form. Note that $\det W\neq 0\neq \det Y$ 
and $\det P'\neq 0\neq \det R'$ due to 
rk$(M_1-I)=$rk$(M_2-I)=2n/3$. Hence, $\det P\neq 0\neq \det R$.~~~~$\Box$

{\bf Proof of Lemma \protect\ref{dimA}:}
Recall that one has $N_2(I+G+H)N_1=M_1M_2$ and that the matrix $M_1M_2$ is 
diagonalizable with  
three eigenvalues each of multiplicity $n/3$. Hence, the quadruple of matrices 
$N_2$, $I+G+H$, $N_1$ and $(M_1M_2)^{-1}$ (their product is 
$I$) is reducible -- if the map $\Psi$ is applied to the quadruple, 
then one obtains a 
quadruple of conjugacy classes of size $2n/3$ the first three of which are 
each with a single 
eigenvalue and with $n/3$ Jordan blocks of size 2 and the fourth of which 
is diagonalizable, with two eigenvalues each of multiplicity $n/3$. 
One can apply the basic technical tool to such a quadruple and deform it into 
one with relatively generic but not generic eigenvalues and in which all four 
matrices are diagonalizable and have two eigenvalues of multiplicity $n/3$. 
This is a quadruple from Case A) (recall that the value of $\xi$ is 
preserved), hence, block-diagonal up 
to conjugacy with diagonal blocks of one and the same size 
(Remark~\ref{caseArem}).

Hence, there exist only block-diagonal up to conjugacy 
quadruples of matrices $N_2$, $I+G+H$, $N_1$ and $(M_1M_2)^{-1}$ and all their 
diagonal blocks are of the same size. The dimension of such a matrix algebra 
is $\leq n^2/2$ with equality if and only if there two diagonal blocks.

\section{Case D)\protect\label{seccaseD}}

Set $s=n/l$ ($l$ was defined in Subsection~\ref{GEPsi}). Hence,
$n=6ks$, $k>1$ and the MVs of $M_1$, $M_2$, $M_3$ equal respectively 
$(sk,sk,sk,sk,sk,sk)$, $(2sk,2sk,2sk)$, $(3sk,3sk)$. 
Case D) can be reduced to Case B)  
like this: if the DSP is solvable in case D), then using Lemma~\ref{helpful} 
one can 
choose the eigenvalues of $M_3$ to be $\pm 1$, i.e. $(M_3)^2=I$, and 
the ones of $M_1$ to form
three couples of opposite eigenvalues; hence, the MV of $(M_1)^2$ is 
$(2sk,2sk,2sk)$ and one has $(M_1)^{-2}=M_2(M_3M_2M_3)$.

Hence, the three 
matrices $(M_1)^2$, 
$M_2$ and $M_3M_2M_3=(M_3)^{-1}M_2M_3$ are from Case B). By assumption, 
they define a block 
diagonal matrix algebra ${\cal A}$ with $2k$ diagonal blocks $3s\times 3s$ 
(Remark~\ref{caseArem}). Hence, dim${\cal A}\leq 18ks^2$. 
The algebra ${\cal A}$ contains the matrices $(M_1)^2$, $(M_3)^2$, 
$(M_1)^{-1}M_3=M_2$ and $M_3(M_1)^{-1}=M_3M_2M_3$. Every matrix from the 
algebra ${\cal B}$ generated by $(M_1)^{-1}$ and $M_3$ (this is also the 
algebra generated by $M_1$, $M_2$ and $M_3$) is representable as 
$K+M_1L+M_3N$, $K,L,N\in {\cal A}$. Hence, 
dim${\cal B}\leq 54ks^2<n^2=36k^2s^2$ and this cannot be $gl(n,{\bf C})$. By 
the Burnside theorem, ${\cal B}$ is reducible.

\section{Proof of Theorem~\protect\ref{diagbasicres} in the case of 
matrices $A_j$\protect\label{matricesA_j}}

Suppose that the Deligne-Simpson problem is weakly solvable in one of cases 
A) -- D) for matrices $A_j$ with relatively generic but not generic 
eigenvalues. By Lemma~\ref{TC=Irr} it is solvable as well. 

Construct a Fuchsian system with matrices-residua from an irreducible triple 
or quadruple corresponding to one of the four cases and with relatively 
generic eigenvalues. One can multiply the 
matrices-residua by $c^*\in {\bf C}$ so that no two eigenvalues differ by a 
non-zero integer and the eigenvalues of the 
monodromy operators become relatively generic. 

Hence, the monodromy group of the system is irreducible. Indeed, 
if it were reducible, then the eigenvalues of the 
diagonal blocks would satisfy only the basic non-genericity relation and its 
corollaries. The sum of the corresponding eigenvalues of the matrices-residua 
is 0 and, hence, one can conjugate simultaneously the matrices-residua to a 
block upper-triangular form, see \cite{Bo2}, Theorem 5.1.2. 
The irreducibility of the monodromy group contradicts 
part 1) of Theorem~\ref{diagbasicres}.

\end{document}